\documentclass[onecolumn,final,a4paper,sort&compress]{elsarticle}

\usepackage[utf8]{inputenc}
\usepackage{graphicx}
\usepackage{wrapfig}
\usepackage{amsfonts}
\usepackage{amsmath}
\numberwithin{equation}{section}
\usepackage{mathtools}
\usepackage{xcolor}
\usepackage{makecell}

\newcommand{\diff}{\mathop{}\!d}
\usepackage[ruled,vlined,linesnumbered]{algorithm2e}
\let\oldnl\nl
\newcommand{\nonl}{\renewcommand{\nl}{\let\nl\oldnl}}
\newcommand{\diffmat}{\mathbf{D}}
\newcommand{\intpmat}{\mathbf{Ip}}

\newcommand{\angmat}{\mathbf{Q}}
\newcommand{\tenprod}{\circ}
\newcommand{\ten}[1]{\mathcal{#1}}
\newcommand{\mat}[1]{\mathbf{#1}}
\renewcommand{\vec}[1]{\mathbf{#1}}
\newcommand{\bytwo}{2\times \cdots \times 2}

\newcommand{\TT}[1]{TT {#1}}
\newcommand{\QTT}[1]{QTT {#1}}
\newcommand{\HAT}[1]{\widetilde{}}

\usepackage{color}
\newcommand{\RED}[1]{\textcolor{red}{#1}}
\newcommand{\NOTE}[1]{\noindent\RED{\textbf{#1}}}

\usepackage[a4paper]{geometry}
\DeclarePairedDelimiter\abs{\lvert}{\rvert}%
\DeclareMathOperator*{\argmin}{argmin} 

\begin{document}

\begin{frontmatter}
  
\title{Tensor Networks for Solving Realistic Time-independent Boltzmann Neutron Transport Equation}

\author[A]{Duc P. Truong}
\author[B]{Mario I. Ortega}
\author[A]{Ismael Boureima}
\author[A]{Gianmarco Manzini}
\author[A]{\\Kim \O. Rasmussen}
\author[A]{and~Boian S. Alexandrov}
  
\address[A]{Theoretical Division, Los Alamos National Laboratory,
  Los Alamos, New Mexico}
  \address[B]{Computer, Computational, and Statistical Sciences Division,
  Los Alamos National Laboratory,
  Los Alamos, New Mexico}
 
  \begin{abstract}
    Tensor network techniques, known for their low-rank approximation
    ability that breaks the curse of dimensionality, are emerging as a
    foundation of new mathematical methods for ultra-fast numerical
    solutions of high-dimensional Partial Differential Equations
    (PDEs).
    Here, we present a mixed Tensor Train (TT)/Quantized Tensor Train
    (QTT) approach for the numerical solution of time-independent
    Boltzmann Neutron Transport equations (BNTEs) in Cartesian
    geometry.
    Discretizing a realistic three-dimensional (3D) BNTE by $(i)$
    diamond differencing, $(ii)$ multigroup-in-energy, and $(iii)$
    discrete ordinate collocation leads to huge generalized eigenvalue
    problems that generally require a matrix-free approach and large
    computer clusters.
    Starting from this discretization, we construct a TT
    representation of the PDE fields and discrete operators, followed
    by a QTT representation of the TT cores and solving the tensorized
    generalized eigenvalue problem in a fixed-point scheme with tensor
    network optimization techniques.
    We validate our approach by applying it to two realistic examples
    of 3D neutron transport problems, currently solved by the PARallel
    TIme-dependent SN (PARTISN) solver\footnote{The "SN method"
    (structured Newton method) is an approximation method for solving
    the radiation transport equation by discretizing both the
    XYZ-domain and the angular variables that specify the direction of
    radiation, developed by Subrahmanyan Chandrasekhar.}.
    We demonstrate that our TT/QTT method, executed on a standard
    desktop computer, leads to a yottabyte compression of the memory
    storage, and more than 7500 times speedup with a discrepancy of
    less than $10^{-5}$ when compared to the PARTISN solution.
  \end{abstract}
 

\end{frontmatter}

\renewcommand\thefootnote{\textcolor{red}{\bf{\{\arabic{footnote}\}}}}

\section{Introduction}
Realistic simulations in classical and quantum physics and chemistry
and complex engineering problems often seek numerical solutions of
high-dimensional partial differential equations (PDEs) or
integro-differential equations.
To numerically solve such equations, discretized versions of the
mathematical models are required.
The discretization represents the solution, a multivariate function,
by its values at a number of grid points and the derivatives of this
function through differences in these values.

The number of grid points increases exponentially with the number
of dimensions, $d$, often making numerical solutions infeasible.
This phenomenon is called the \emph{curse of dimensionality}
\cite{bellman1966dynamic}, it causes the poor computational scaling of
numerical algorithms, and is the primary challenge for
multidimensional numerical computations regardless of the specific
problem.
Importantly, exascale high-performance computing, which offers
strategies for optimization, cannot break the curse of dimensionality.

A very recent and promising approach to mitigate or even remove the
curse of dimensionality is based on \emph{Tensor Networks (TNs)},
which represent a non-trivial generalization of tensor factorization.
In \emph{big data analytics}, TNs allow the examination of
high-dimensional data by partitioning them into smaller, manageable
blocks, i.e., by approximating a high-dimensional array by a network
of low-dimensional tensors.
In this sense, we can consider classical tensor factorization formats
such as \emph{Canonical Polyadic Decomposition (CPD)} and \emph{Tucker
decomposition} as TNs.
Tensor networks (TNs), which are multilinear-algebra data structures
\cite{cichocki2016tensor,bachmayr2016tensor}, have been developed and enjoyed great
success in theoretical physics
\cite{penrose1971applications,feynman2023quantum,deutsch1989quantum},
as well as in data science in machine learning
\cite{cichocki2017tensor,cichocki2014tensor}.
TNs are emerging as a promising novel strategy for breaking the curse
of dimensionality \cite{hackbusch2012tensor} in the numerical solution
of high-dimensional differential equations and integrals
\cite{khoromskij2018tensor,alexandrov2023challenging,gourianov2022quantum}.
Several authors have recently used TNs seeking fast and accurate
solutions of specific examples of equations, such as the Poisson
equation~\cite{chipot2022numerical}, the Schrödinger
equation~\cite{gelss2022solving}, the Poisson-Boltzmann
equation~\cite{kweyu2021solution}, the Smoluchowski
equation~\cite{matveev2016tensor, manzini2021nonnegative}, the Maxwell
equations~\cite{manzini2023tensor}, the Vlasov–Poisson
equations~\cite{kormann2015semi}, the neutron diffusion equation
\cite{kusch2022low}, and others.

Numerical algorithms for solving PDEs are, in general, composed of two
main ingredients.
The first one is the \emph{grid functions}, which are the value of the
functions evaluated at the nodes of a discrete grid.
The second ingredient is given by the \emph{discrete operators}, which
are the discrete analogs of differential operators (e.g., gradient,
curl, divergence), a combination of them (e.g., Laplacian, div-grad,
etc), and numerical interpolation and integration operators.

The Boltzmann Neutron Transport Equation (BNTE)
\cite{lewis1984computational} is an integro-differential equation for
the neutron angular flux describing the location in space, energy, and
the direction of neutrons in a physical system of interest, such as a
nuclear reactor.
The BNTE discretization is driven by the physics of neutron
interactions and generates extremely large linear
systems of equations that cannot be straightforwardly handled
numerically.
Here, we introduce a new approach to solving the BNTE, which is based
on the Tensor Train (TT) format~\cite{oseledets2011tensor} combined
with the Quantized Tensor Train (QTT) format~\cite{khoromskij2011d},
and applied to the finite difference discretization method.
This discretization approach leads to ultra-large linear systems of equations 
and generalized eigenvalue problems.
\TT{solvers} exists to solve linear systems in \TT{format}.
Examples are the \emph{Alternating Linear Scheme
(ALS)}~\cite{holtz2012alternating}, the \emph{Density Matrix
Renormalization Group Algorithm (DMRG)}~\cite{savostyanov2011fast},
the \emph{Alternating Minimal Energy (AMEn)
methods}~\cite{dolgov2014alternating}, etc.
In contrast, methods for solving eigenvalue problems in TT/QTT format are
still an open research topic with a few algorithms being
proposed~\cite{dolgov2014computation,ruymbeek2022subspace}.
In this work, we do not utilize a linear solver based on Krylov
subspaces but instead, we seek the solution directly in the
\TT{format} \cite{oseledets2012solution}, while we redesign a
fixed-point scheme in the \QTT{format} to compute the largest
eigenvalue of the BNTE problem.
We compare the efficiency, speed, memory requirements, and accuracy,
of our method to the traditional matrix-free approach
\cite{wang2015matrix} for solving realistic BNTEs for three-dimensional
systems \cite{ortega2020rayleigh}.
We demonstrate that our TT/QTT method can give us a more than 7500
times speedup, based on yottabyte compression, while preserving an
accuracy of $10^{-5}$.

The outline of the paper is as follows.
In Section~\ref{sec2:Boltzmann-NTE}, we review some basic concepts and
the matrix formulation of the Boltzmann Neutron Transport Equation
that leads to a generalized eigenvalue problem.
In Section~\ref{sec3:Tensor-Networks}, we introduce the tensor
notation and the definitions of tensor networks, TT and QTT formats,
and their application to PDEs.
In Section~\ref{sec4:tensorization-Boltzmann-NTE}, we present our
TT/QTT design of the numerical solution of the Boltzmann Neutron
Transport Equation and introduce our fixed-point algorithms in tensor train format.

In Section~\ref{sec5:numerical-experiments}, we present our numerical
results and assess the performance of our method, comparing its
efficiency to the efficiency of the PARallel TIme Dependent SN
(PARTISN) solver \cite{alcouffe2005partisn} applied to the same
problems.
In Section~\ref{sec:conclusions}, we offer our final remarks and
discusses possible future work.

For completeness, we report part of the details the BNTE
matrix representation in \ref{eqn:Appendix} and details of the tensorization approach in \ref{Appendix_B}.

\section{Methods}
\subsection{Boltzmann Neutron Transport Equation}
\label{sec2:Boltzmann-NTE}

In the design of nuclear systems such as nuclear reactors, nuclear
engineers are especially interested in determining the criticality
conditions.
\emph{Criticality} is the ability of a nuclear system to sustain a
nuclear chain reaction (number of neutrons created in fission is equal
to the number of neutrons lost in the system) without an external
source \cite{bell_nuclear_1970}.
To determine the criticality of a system, we consider the
time-independent, \emph{k-effective eigenvalue problem} and the
\emph{alpha-effective eigenvalue problem} for the eigenpairs
$\big(k_{\text{eff}},\psi\big)$ and $\big(\alpha,\psi\big)$.
Both eigenvalue problems provide helpful insight into nuclear systems
and are used throughout nuclear reactor design and dynamics,
criticality safety, and nuclear non-proliferation applications.
In both cases, the eigenfunction $\psi$ is a function of the position
vector $\vec{r}$, the direction variable is $\hat{\Omega}$, and the
group energy $E$, i.e., $\psi(\vec{r},\hat{\Omega},E)$.
The position vector $\vec{r}$ varies on the space domain
\begin{align*}
  \mathcal{D}\equiv
  \Big\{
  \vec{r} = (x,y,z) \in \mathbb{R}^3 |
  a_{x}\leq x\leq b_{x},
  a_{y}\leq y\leq b_{y},
  a_{z}\leq z\leq b_{z}
  \Big\},
\end{align*}
where $\big[a_{x},b_{x}\big]$,, $\big[a_{y},b_{y}\big]$,,
$\big[a_{z},b_{z}\big]$ are bounded subintervals of $\mathbb{R}$;
the direction variable $\hat{\Omega}$ varies on $\mathcal{S}^{2}$, the
unit sphere in $\mathbb{R}^{3}$;
the energy group varies, for convenience, between a minimum and
maximum values, e.g., $E\in\big[E_{\text{min}},E_{\text{max}}\big]$.
For both problems, we assume vacuum Dirichlet boundary conditions:
\begin{align}
  \psi( \vec{r}, \hat{\Omega}, E ) = 0
  \text{ for all } \vec{r} \in \partial \mathcal{D}
  \text{ and } \hat{\Omega} \in \mathcal{S}^{2}
  \text{ with } \vec{n}(\vec{r}) \cdot \hat{\Omega} < 0,
  \label{eq:BCS}
\end{align}
where $\vec{n}(\vec{r})$ is the unit normal vector to
$\partial\mathcal{D}$, the boundary of the computational domain
$\mathcal{D}$, and pointing out of $\mathcal{D}$

In both eigenvalue problems, the physics of neutron interactions with
matter is captured through \emph{neutron cross sections},
probabilistic measures of certain nuclear reactions taking
place~\cite{duderstadt_nuclear_1976}.
The cross section are denoted as $\sigma$, $\sigma_{s}$, $\sigma_{f}$,
the total, scattering, and fission cross sections, respectively, and
are functions of material, energy, and collision angle.
We also use the symbol $\nu$ 
to denote the number of
neutrons emitted in fission, while $\chi(E'\rightarrow E )$ denotes the
probability distribution function expressing the probability that a
neutron with energy $E'$ induces fission on fissile nuclei and creates
neutrons with energy $E$.

\subsubsection{The k-effective eigenvalue problem}
For the eigenpair $\big(k_{\text{eff}},\psi\big)$, we consider the
eigenvalue problem
\begin{align}
  &\hat{\Omega}\cdot\nabla\psi(\vec{r},\hat{\Omega},E)
  + \sigma(\vec{r},E)\psi( \vec{r},\hat{\Omega},E)
  \nonumber\\[0.5em]
  &\qquad\qquad
  =\int_{\mathcal{S}^{2}}\diff\hat{\Omega'}\,\psi(\vec{r},\hat{\Omega}',E)\sigma_{s}(\vec{r},\hat{\Omega}\cdot\hat{\Omega}',E)
  \nonumber\\[0.5em]
  &\qquad\qquad\quad
  +\frac{1}{k_{\text{eff}}}\int_{0}^{\infty}\diff E'\int_{\mathcal{S}^{2}}\diff\hat{\Omega'}\chi(E'\rightarrow E)\nu\sigma_{f}(\vec{r},E')\psi(\vec{r},\hat{\Omega}',E').
  \label{eq:keff3DNeutronTransport}
\end{align}
In~\eqref{eq:keff3DNeutronTransport}, the eigenvalue $k_{\text{eff}}$
scales the number of neutrons emitted in fission, and its value
determines the criticality of the system:
\begin{align*}
  k_{\text{eff}}\quad
  \begin{cases}
    \,\, > 1, & \text{supercritical,}\\[0.25em]
    \,\, = 1, & \text{critical,}\\[0.25em]
    \,\, < 1, & \text{subcritical.}
  \end{cases}
\end{align*}
A system with $k_{\text{eff}} = 1$ is considered critical and the
number of neutrons is constant in time.
The number of neutrons goes to zero or infinity for subcritical or
supercritical systems, respectively.

\subsubsection{The alpha-effective eigenvalue problem}
For the eigenpair $\big(\alpha,\psi\big)$, we consider the
eigenvalue problem
\begin{align}
  & \frac{\alpha}{v(E)} \psi( \vec{r}, \hat{\Omega}, E ) + \hat{\Omega} \cdot \nabla\psi( \vec{r}, \hat{\Omega}, E ) + \sigma(\vec{r},E) \psi( \vec{r}, \hat{\Omega}, E )
  \nonumber\\[0.5em]
  &\qquad\qquad
  = \int_{\mathcal{S}^{2}} \diff \hat{\Omega'} \, \psi( \vec{r}, \hat{\Omega}', E ) \sigma_{s}(\vec{r},\hat{\Omega} \cdot \hat{\Omega}', E)
  \nonumber\\[0.5em]
  &\qquad\qquad\quad
  + \int_{0}^{\infty} \diff E' \int_{\mathcal{S}^{2}} \diff \hat{\Omega'} \chi(E' \rightarrow E ) \nu \sigma_{f}(\vec{r}, E') \psi( \vec{r}, \hat{\Omega}', E' ).
  \label{eq:alpha3DNeutronTransport}
\end{align}
The alpha-eigenvalue $\alpha$ gives a measure of the asymptotic time behavior
of a system and can also be used to determine the criticality of a
system:
\begin{align*}
  \alpha\quad
  \begin{cases}
    \,\, > 0, & \text{supercritical,}\\[0.25em]
    \,\, = 0, & \text{critical,}\\[0.25em]
    \,\, < 0, & \text{subcritical.}
  \end{cases}
\end{align*}
For supercritical systems, the alpha-eigenvalue is the $e$-folding
time for the neutron angular flux in a nuclear system, i.e., the
  time interval in which $\psi$ grows exponentially by a factor of
  e.

\subsubsection{Discretization of the Three-Dimensional Neutron Transport}

The numerical solution of equations \eqref{eq:keff3DNeutronTransport}
and \eqref{eq:alpha3DNeutronTransport} requires discretization along
the spatial, angular, and energy variables.
In this section, we describe the major discretization issues and the
matrix formulation, which is similar to that of
Refs. \cite{brown_linear_1995} and \cite{ortega_2020_alpharq}.
Specifically, we first consider the energy variable; then, the angular
variable; and, finally, the space variables.

\subsubsection{Discrete formulation of k-effective and alpha-effective eigenvalue problems}

The discretized forms of equations \eqref{eq:keff3DNeutronTransport}
and \eqref{eq:alpha3DNeutronTransport} are (see \ref{eqn:Appendix} for details):
\begin{multline}
  \frac{\mu_{\ell}}{4 \Delta x_{i}} \Bigg [ \sum_{k'=k-1}^{k} \sum_{j'=j-1}^{j} \psi_{g,\ell,k',j',i} - \psi_{g,\ell,k',j',i-1} \Bigg ] + \\
  \frac{\eta_{\ell}}{4 \Delta y_{j}} \Bigg [ \sum_{k'=k-1}^{k} \sum_{i'=i-1}^{i} \psi_{g,\ell,k',j,i'} - \psi_{g,\ell,k',j-1,i'} \Bigg ] + \\
  \frac{\xi_{\ell}}{4 \Delta z_{k}} \Bigg [ \sum_{j'=j-1}^{j} \sum_{i'=i-1}^{i} \psi_{g,\ell,k,j',i'} - \psi_{g,\ell,k-1,j',i'} \Bigg ] + \\
  \frac{\sigma_{g,k,j,i}}{8} \Bigg [ \sum_{k'=k-1}^{k} \sum_{j'=j-1}^{j} \sum_{i'=i-1}^{i} \psi_{g,\ell,k',j',i'} \Bigg ] = \\
  \frac{1}{k_{\text{eff}}}\frac{1}{8} \sum_{g'=1}^{G} \chi_{g'g}\nu\sigma_{f,g',k,j,i} \sum_{\ell' = 1}^{L} w_{\ell'} \Bigg [ \sum_{k'=k-1}^{k} \sum_{j'=j-1}^{j} \sum_{i'=i-1}^{i} \psi_{g',\ell',k',j',i'} \Bigg ] + \\
  \frac{1}{8} \sum_{g'=1}^{G} \sigma_{s,g,g',k,j,i} \sum_{\ell' = 1}^{L} w_{\ell'} \Bigg [ \sum_{k'=k-1}^{k} \sum_{j'=j-1}^{j} \sum_{i'=i-1}^{i} \psi_{g',\ell',k',j',i'} \Bigg ],
  \label{eq:Discretized_kNTE}
\end{multline}
and
\begin{multline}
  \frac{1}{8}\frac{\alpha}{v_{g}} \Bigg [ \sum_{k'=k-1}^{k} \sum_{j'=j-1}^{j} \sum_{i'=i-1}^{i} \psi_{g,\ell,k',j',i'} \Bigg ] + \\ \frac{\mu_{\ell}}{4 \Delta x_{i}} \Bigg [ \sum_{k'=k-1}^{k} \sum_{j'=j-1}^{j} \psi_{g,\ell,k',j',i} - \psi_{g,\ell,k',j',i-1} \Bigg ] + \\
  \frac{\eta_{\ell}}{4 \Delta y_{j}} \Bigg [ \sum_{k'=k-1}^{k} \sum_{i'=i-1}^{i} \psi_{g,\ell,k',j,i'} - \psi_{g,\ell,k',j-1,i'} \Bigg ] + \\
  \frac{\xi_{\ell}}{4 \Delta z_{k}} \Bigg [ \sum_{j'=j-1}^{j} \sum_{i'=i-1}^{i} \psi_{g,\ell,k,j',i'} - \psi_{g,\ell,k-1,j',i'} \Bigg ] + \\
  \frac{\sigma_{g,k,j,i}}{8} \Bigg [ \sum_{k'=k-1}^{k} \sum_{j'=j-1}^{j} \sum_{i'=i-1}^{i} \psi_{g,\ell,k',j',i'} \Bigg ] = \\
  \frac{1}{8} \sum_{g'=1}^{G} \chi_{g'g}\nu\sigma_{f,g',k,j,i} \sum_{\ell' = 1}^{L} w_{\ell'} \Bigg [ \sum_{k'=k-1}^{k} \sum_{j'=j-1}^{j} \sum_{i'=i-1}^{i} \psi_{g',\ell',k',j',i'} \Bigg ] + \\
  \frac{1}{8} \sum_{g'=1}^{G} \sigma_{s,g,g',k,j,i} \sum_{\ell' = 1}^{L} w_{\ell'} \Bigg [ \sum_{k'=k-1}^{k} \sum_{j'=j-1}^{j} \sum_{i'=i-1}^{i} \psi_{g',\ell',k',j',i'} \Bigg ],
  \label{eq:Discretized_alphaNTE}
\end{multline}
for
$g = 1, \dots, G$,
$i = 1, \dots, M$,
$j = 1, \dots, J$,
$k = 1, \dots, K$, and
$\ell = 1, \dots, L$.
In equations \eqref{eq:Discretized_kNTE} and
\eqref{eq:Discretized_alphaNTE}, $\sigma_{g}$, $\sigma_{s,g,g'}$,
$\nu\sigma_{f,g}$, and $v_{g}$, the total, scattering, and fission
cross sections and the velocity for energy group $g$, are assumed to
be constant on the cell $(x_{i-1} < x < x_{i}, y_{j-1} < y < y_{j},
z_{k-1} < z < z_{k})$ (see~\ref{eqn:Appendix}).
We denote these values on each cell by $\sigma_{g,k,j,i}$,
$\sigma_{s,g,g',k,j,i}$, $\nu\sigma_{f,g,k,j,i}$, and $v_{g}$.

\subsubsection{Matrix formulation and the fixed-point iterative scheme}

We use a matrix formulation similar to that of
Refs. \cite{brown_linear_1995} and \cite{ortega_2020_alpharq} to rewrite equations \eqref{eq:Discretized_kNTE} and
\eqref{eq:Discretized_alphaNTE} in a more compact way:
\begin{equation}
\mat{H} \Psi = \bigg [ \mat{S} + \frac{1}{k_{\text{eff}}} \mat{F} \bigg ] \Psi,
\label{eq:OpFormkNTE}
\end{equation}
and
\begin{equation}
\bigg [ \alpha \mat{V}^{-1} + \mat{H} \bigg ] \Psi = \bigg [ \mat{S} + \mat{F} \bigg ] \Psi.
\label{eq:OpFormAlphaNTE}
\end{equation}
Equations \eqref{eq:OpFormkNTE} and \eqref{eq:OpFormAlphaNTE} are
linear equations with $G\times L\times K\times J\times M$ unknowns for
$\Psi$.
The matrix operators, $\mat{H}$, $\mat{S}$, $\mat{F}$, and $\mat{V}$
are constructed as Kronecker products of various smaller matrices that
couple spatial, angular, and energy angular flux unknowns.

Starting from random values for $\Psi^{(0)}$ and $k^{(1)}$,
we solve
the k-effective problem~\eqref{eq:OpFormkNTE} through the iterative
process for $\tau\geq1$:
\begin{equation}
  \begin{aligned}
    \Psi^{(\tau+1)} &= \mat{H}^{-1}\bigg [ \mat{S} + \frac{1}{k_{\text{eff}}^{(\tau)}} \mat{F} \bigg ] \Psi^{(\tau)},\\[0.5em]
    k_{\text{eff}}^{(\tau+1)} &= k_{\text{eff}}^{(\tau)} \frac{ \sum \mat{F} \Psi^{(\tau+1)}}{ \sum \mat{F} \Psi^{(\tau)}}.
  \end{aligned}
  \label{eqn:k_eff_fixed_points}
\end{equation}

Starting from random values for $\Psi^{(0)}$ and setting
$\alpha^{(0)}=0$, we solve for $k_{\text{eff}}^{(0)}$. We then set $\alpha^{(1)}=0.01$ and solve for $k_{\text{eff}}^{(1)}$. Continuing this process we solve
problem~\eqref{eq:OpFormAlphaNTE} through the iterative process for
$\tau\ge1$:
\begin{equation}
  \begin{aligned}
    \bigg [ \alpha^{(0)}\mat{V}^{-1} + \mat{H} \bigg ]\Psi^{(0)} &= \bigg [ \mat{S} +  \frac{1}{k_{\text{eff}}^{(0)}} \mat{F} \bigg ] \Psi^{(0)} \leftarrow \text{solve this k-effective with } \alpha^{(0)},\\
    \bigg [ \alpha^{(1)}\mat{V}^{-1} + \mat{H} \bigg ]\Psi^{(1)} &= \bigg [ \mat{S} +  \frac{1}{k_{\text{eff}}^{(1)}} \mat{F} \bigg ] \Psi^{(1)} \leftarrow \text{solve this k-effective with } \alpha^{(1)},\\
    \alpha^{(\tau+1)} &= \alpha^{(\tau)} + \dfrac{1-{k_{\text{eff}}^{(\tau-1)}}}{(k_{\text{eff}}^{(\tau)}-k_{\text{eff}}^{(\tau-1)})*(\alpha^{(\tau)} - \alpha^{(\tau-1)})},\\
    \bigg [ \alpha^{(\tau+1)}\mat{V}^{-1} + \mat{H} \bigg ]\Psi^{(\tau+1)} &= \bigg [ \mat{S} +  \frac{1}{k_{\text{eff}}^{(\tau+1)}} \mat{F} \bigg ] \Psi^{(\tau+1)} \leftarrow \text{solve this k-effective with } \alpha^{(\tau+1)}.
  \end{aligned}
  \label{eqn:alpha_eig_fixed_points}
\end{equation}
Since the matrices, $\mat{H}$, $\mat{S}$, $\mat{F}$, and $\mat{V}$,
tend to be very large; in most cases, this iterative process is
accomplished in PARTISN using a matrix-free method such as those
described in Ref.~\cite{lewis_computational_1984}.



\subsection{Tensors and Tensor Networks}
\label{sec3:Tensor-Networks}



\medskip
In our TT/QTT reformulation of the algorithms of
Section~\ref{sec2:Boltzmann-NTE}, we will make extensive use of real,
multidimensional tensors, i.e., multidimensional arrays of real
numbers.
We refer the reader to Refs.~\cite{kolda2009tensor,oseledets2011tensor} and \ref{Appendix_B} for
more details about the notation and the concepts that we briefly
review in this section.

\subsubsection{Tensor network formats and tensor factorizations}
\label{subsec:TNs}


The total number of elements $N$ of a $d$-dimensional tensor, with
$n_k=\mathcal{O}(n)$ elements in each dimension $k=1,2,\ldots d$ is
exponential in $d$, i.e., $N=\mathcal{O}(n^d)$.
Approximate tensor factorization compresses the full tensor with
``acceptable'' accuracy by using much fewer elements.
Such a factorization is achieved through a multidimensional
minimization that can include various constraints (sparsity,
non-negativity, etc.)~\cite{kolda2009tensor}.
To define the tensor decompositions, we need the notion of
\emph{tensor rank}.
A $d$-dimensional, rank-1 tensor is a tensor that can be represented as tensor product of $d$ vectors,
e.g., $\ten{G}=\vec{g}_1\circ\vec{g}_2\circ\ldots\circ\vec{g}_d$, or,
componentwise,
$\ten{G}(i_1,i_2,\ldots,i_d)=\vec{g}_1(i_1)\vec{g}_2(i_2)\ldots\vec{g}_d(i_d)$.
The canonical rank of a tensor is the minimal number $R$ of rank-1
tensors whose sum is equal to this tensor.
\begin{figure}[h] 
  \centering
  \includegraphics[width=0.75\textwidth]{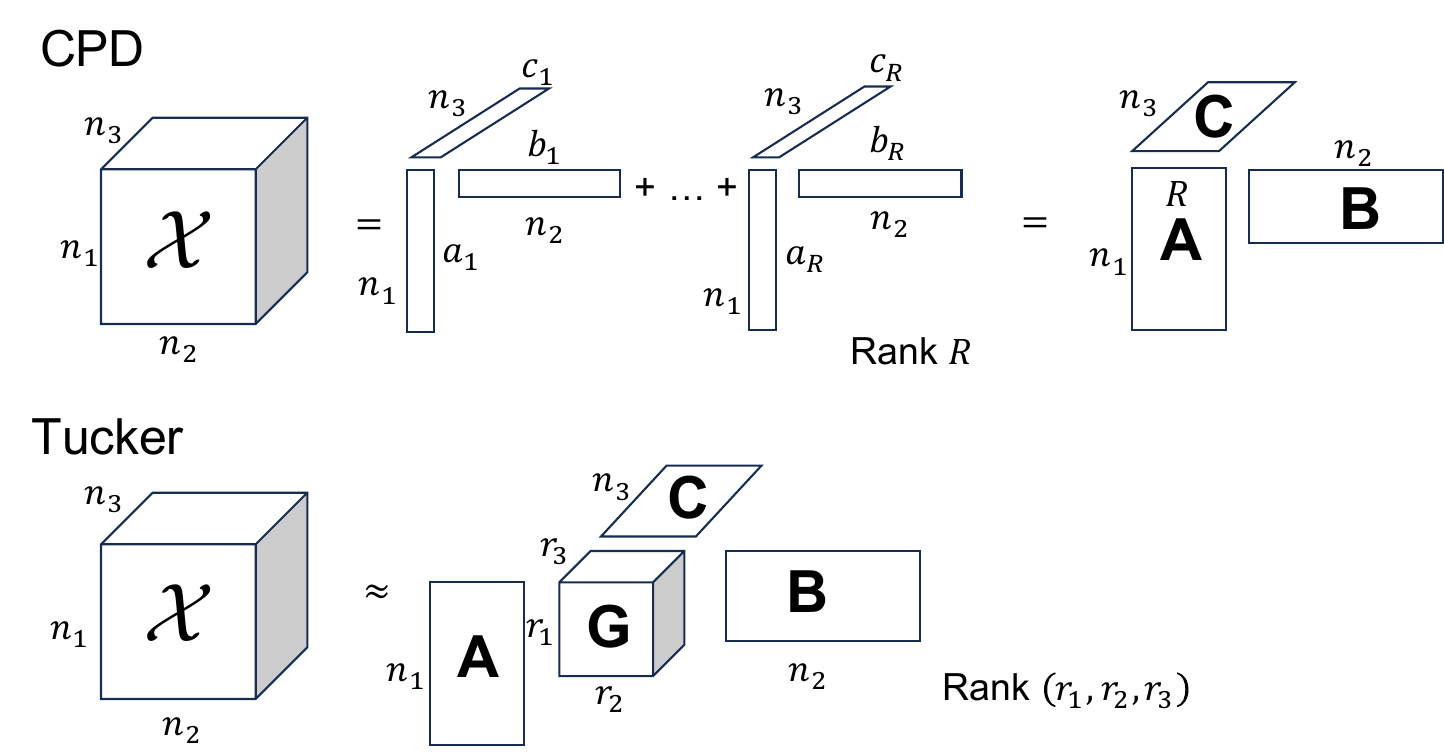}
  \caption{\textit{Top panel:} Canonical Polyadic Decomposition (CPD)
    of a 3D tensor with size $N_{\chi}=n_1 \times n_2 \times n_3$ and
    rank $R$. \textit{Bottom panel:} Tucker decomposition of a 3D
    tensor with size $N_{\chi}=n_1\times n_2\times n_3$ and multi-rank
    rank $r =[r_1,r_2,r_3]$ .}
    \vspace{4mm}
    \label{fig:cpd_tucker}
\end{figure}

As shown in Fig.~\ref{fig:cpd_tucker} (top panel), the Canonical
Polyadic Decomposition factorizes a $d$-dimensional tensor with rank
$R$, presenting it as a sum of $R$ rank-1, $d$-dimensional tensors,
cf.~\cite{harshman1972determination}.
This decomposition has the smallest possible number of elements,
$\ten{O}(NdR)$.
However, it requires knowledge of the canonical rank, 
whose computation is an NP-hard problem \cite{haastad1989tensor}.
Therefore, the approximation of the full tensor by CPD, with its canonical rank $R$, can be inaccurate or ill-posed \cite{de2008tensor}.
As shown in~Fig.~\ref{fig:cpd_tucker} (bottom panel), \emph{Tucker
decomposition} factorizes a $d$-dimensional tensor by a product of $d$
factor matrices and a small core tensor, $\ten{G}$ with dimensions
$r_1\times r_2\times\ldots\times r_d$, see
Refs. \cite{de2000multilinear,tucker1966some}.
The set of the dimensions of the core tensor, $\ten{G}$,
$\mathbf{r}:=\big[r_1,r_2,\dots,r_d\big]$, is called \emph{Tucker
multi-rank}.
The number of Tucker decomposition elements, $\ten{O}(Ndr +r^d)$,
remains exponential in $d$.

The Tensor Train (TT) format~\cite{oseledets2011tensor}, seen as a
linear chain of products of 3D tensors, is a very effective
alternative to CPD and Tucker decomposition.
Precisely, the \TT{approximation} $\ten{X}^{TT}$ of a
$d$-dimensional tensor $\ten{X}$ is a tensor with elements
\begin{align}
  \ten{X}^{TT}(i_1,i_2,\dots,i_d)
  =  \sum^{r_1}_{\alpha_{1}=1}\ldots\sum^{r_d-1}_{\alpha_{d-1}=1}
  \ten{G}_1(1,i_1,\alpha_1)\ten{G}_2(\alpha_1,i_2,\alpha_2)\ldots\ten{G}_d(\alpha_{d-1},i_d,1)
  + \varepsilon,
  \label{eqn:TT_def_element}
\end{align}
where the last term, $\varepsilon$, is a tensor with the same
dimensions of $\ten{X}$ representing the approximation error.
Equivalently, we can also denote the \TT{format} by the multiple
matrix product

\begin{align}
  \ten{X}^{TT}(i_1,i_2,\dots,i_d) =
  \mat{G}_1(i_1)\mat{G}_2(i_2)\dots\mat{G}_d(i_d) + \varepsilon,
  \label{eqn:TT-vector}
\end{align}
where each term
$\big(\mat{G}_{k}(i_k)\big)_{\alpha_{k-1},\alpha_{k}}$,
$i_k=1,2,\ldots,n_k$, $k=1,2,\ldots,d$, is a matrix of size
$r_{k-1}\times r_k$ (with the assumption that $r_0=r_d=1$).
Therefore, the \TT{cores} $\ten{G}_k(:,i_k,:)$ are a set of matrix
slices $\mat{G}_k(i_k)$ that are labeled with the single index $i_k$.
The entries of the integer array
$\mathbf{r}=\big[r_1,\dots,r_{d-1}\big]$ are the \TT{ranks}, and
quantify the compression effectiveness.
Since each \TT{core} only depends on a single mode index of the full
tensor $\ten{X}$, e.g., $i_k$, the TT format effectively embodies a
discrete separation of variables \cite{bachmayr2016tensor}.
When the \TT{ranks} are relatively small with respect to the problem
size, a TT-based approach is referred to as a \emph{low-rank
approximation} \cite{bachmayr2023low}.

Assuming that $n_k=\mathcal{O}(n)$ and $r_k=\mathcal{O}(r)$ for some
nonnegative integers $n$ and $r$, and for all $k=1,2,\ldots,d$, the
total number of elements that \TT{format} stores is proportional to
$\ten{O}(2nr+(d-2)nr^2)$, which is linear with the number of
dimensions $d$.
\begin{figure}[h] 
  \centering
  \includegraphics[width=0.75\textwidth]{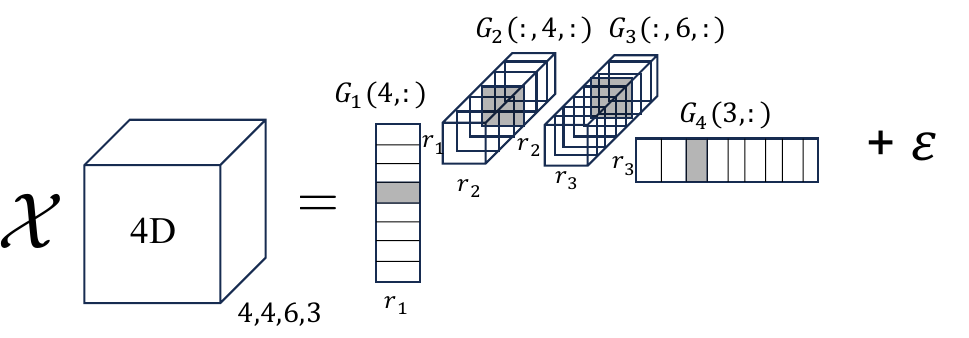}
  \caption{Approximate \TT{decomposition} of a 4D tensor $\ten{X}$,
    with \TT{ranks} $\mathbf{r} = \big[r_1,r_2,r_3,r_4\big]$, and
    approximation error $\varepsilon$, in accordance with
    Eq.~\eqref{eqn:TT-vector}.
    We compute the elements of full tensor $\ten{X}$ as products of
    matrix slices and row and column vectors of the \TT{cores} of
    $\ten{X}^{TT}$.
    For example, element $\ten{X}(4,4,6,3)$ is restored by the product
    of the fourth row vector from the first core $\ten{G}_1$, the
    fourth and sixth matrix slices from the intermediate cores
    $\ten{G}_2$ and $\ten{G}_3$, and the third column vector from the
    fourth and final core $\ten{G}_4$.
    Column and row vectors and matrix slices are grey-shaded in the
    figure.}
    \label{fig:TT_4D}
\end{figure}
Fig.~\ref{fig:TT_4D} illustrates how using the \TT{format}, 
we can approximate a four-dimensional array with a certain error
$\varepsilon$ that we can suitably control.
In fact, for a tensor admitting CPD with rank $R$ and error $\eta$,
there exists an approximate \TT{format} factorization with ranks
$r_k<R$ for each $k=1,2,\ldots,d$ and $\varepsilon<\eta
\sqrt{(d-1)}$ \cite{oseledets2010tt}.

\begin{figure}[t] 
  \centering
  \includegraphics[width=0.75\textwidth]{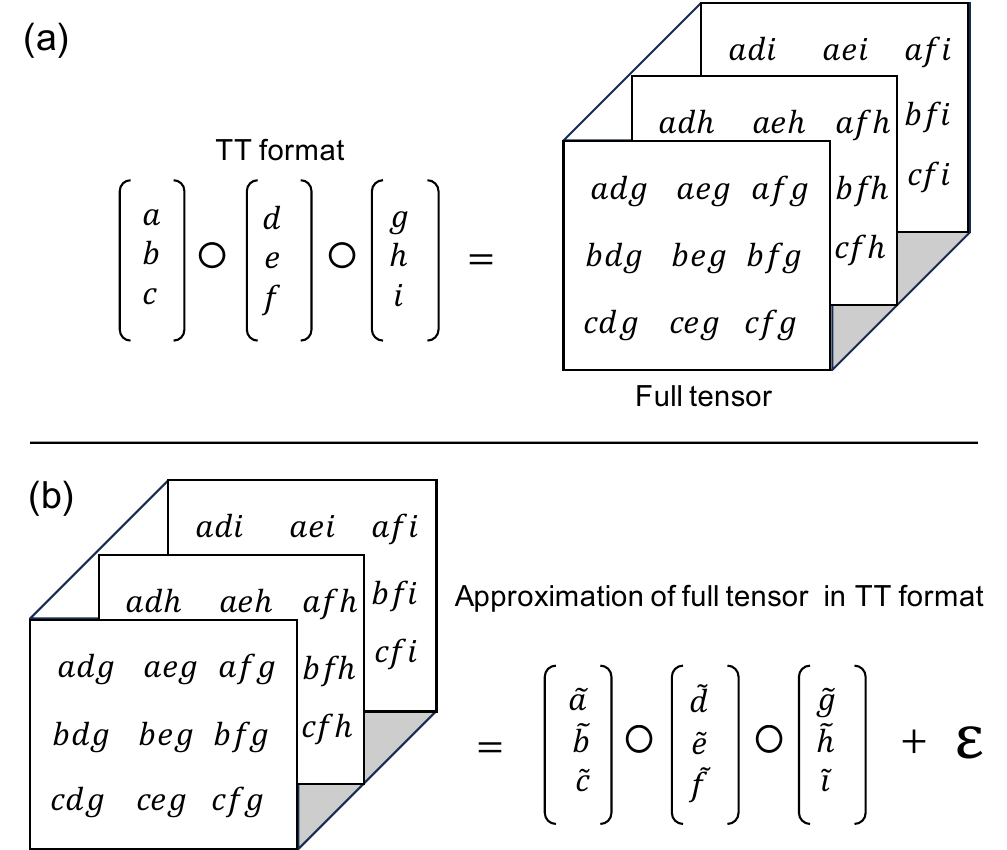}
  \caption{Breaking the curse of dimensionality: $(A)$ The tensor product of three rank-1 vectors is a $3\times3\times3$ tensor, while $(b)$ a
     $3\times3\times3$ tensor can be factorized into a rank-1 product of three vectors 
    with a controlled error $\varepsilon$.}
  \label{fig:structured_unstructured}
\end{figure}

\subsubsection{Grid functions in \TT{format} and curse of dimensionality}
How does the \TT{format} break the curse of dimensionality when
numerically solving high-dimensional PDEs?
The numerical integration of a PDE requires representing 
grid functions and discrete operators in \TT{format}.
In Figure~\ref{fig:structured_unstructured}, a discretization of a 3D
function with three points per dimension, resulting in a dense tensor
(Fig.~\ref{fig:structured_unstructured}a-right), is also presented in
a tensor network format
(Fig.~\ref{fig:structured_unstructured}a-left).
The depicted \TT{format} equals a CPD format with rank $R=1$.
The tensor grid structure makes it possible to store the full-grid
tensor with $n=3$ entries per direction using only $nd = 9$ numbers.
When expanded (Fig.~\ref{fig:structured_unstructured}a-right), the same
amount of information requires storing $n^d=27$ numbers.
In the former situation, the complexity is proportional to $d$ (linear
in $d$), while, in the latter case, it is exponential in $d$.
Therefore, we break the curse of dimensionality if we can perform all
operations of the numerical algorithm solving the PDE through a tensor
network format.
Unfortunately, PDEs are usually discretized using the full tensor
format, (Fig.~\ref{fig:structured_unstructured}b-left), which means
that to break the curse of dimensionality, we must reformulate the PDE
operators, functions, and algorithms on tensor grid functions using
the approximate \TT{format}
(Fig.~\ref{fig:structured_unstructured}b-right) with a controlled error
$\varepsilon$.

\subsubsection{Differentiation, integration and interpolation operators of grid function in \TT{format}}
\label{subsection:TT_operators}
A critical aspect that enables the implementation of the \TT{format}
in numerical PDE methods is the ability to represent all discrete
operators, such as differentiation, interpolation, integration,
multiplication, etc. in \TT{format}
\cite{manzini2023tensor,alexandrov2023challenging}.
For illustration, we consider $d=4$, and a real-valued,
four-dimensional function $f(x_1,x_2,x_3,x_4)$, where each independent
variable $x_k$, $k=1,2,3,4$, is defined over a proper domain range,
e.g., $\mathcal{D}_k$, a bounded subinterval of $\mathbb{R}$.
We introduce a four-dimensional, regular, Cartesian grid covering
domain
$\mathcal{D}=\mathcal{D}_1\times\mathcal{D}_2\times\mathcal{D}_3\times\mathcal{D}_4$
and having $n_k$ nodes along the $k$-th direction.
We let $\ten{F}=\big(\ten{F}(i_1,i_2,i_3,i_4)\big)$ denote the tensor
that consists of the values of $f$ sampled at the grid nodes indexed by
$(i_1,i_2,i_3,i_4)$.
We also let the tensor $\ten{F}$ be approximated by tensor $\ten{G}^{TT}$ in
\TT{format} with cores $\ten{G}_{k}(:,i_k:)$ as
in~\eqref{eqn:TT_def_element} and approximation error $\varepsilon$.
The tensor $\varepsilon$ has the same dimensions and size as $\ten{F}$ and $\ten{G}$ 
and may include other approximation errors due to the numerical differentiation, 
integration, and interpolation and its value may be different at any instance. More details on \TT{format} representation of differentiation, integration, and interpolation operators are given in  \ref{APP:Diff_Operators},
\ref{APP:Int_Operators}, and \ref{APP:Inter_Operators}, respectively.

\subsubsection{Linear Operators in \TT{format}.}
\label{SUB:Lin_operinTT}
A numerical discretization often transforms the unknown multivariate
function into a very long vector of degrees of freedom and the
operators acting linearly on this vector become very large, sparse
matrices.
A linear operator $\mat{A}$ acting on a vector $\mat{x}$ with size
$N_{\mat{x}}=n_1n_2\dots n_d$, which collects the degrees of freedom
of a PDE associated with the multidimensional grid nodes indexed by
$(i_1,i_2,\ldots,i_d)$, has a matrix representation with size
$N_{\mat{A}}=N_{\mat{x}}\times N_{\mat{x}}$.
Instead of using a plain matrix-vector format, we retain the
\TT{format} for matrices representing linear operators (see,
e.g.,~\cite[Section 4.3]{oseledets2011tensor}
and~\cite{oseledets2012solution}),
which generalizes the \TT{format} of~\eqref{eqn:TT_def_element} from
``multi-dimensional vectors'' to ``multi-dimensional matrices''.
To this end, we first note that we can represent the degrees of
freedom of a $d$-dimensional PDE by reshaping vector $\mat{x}$ into a
$d$-dimensional tensor, e.g., $\ten{X}(j_1,\dots,j_d)$.
Consistently, we can reshape operator $\mat{A}$ into a ``matrix''
tensor, e.g., $\ten{A}\big( (i_1,\ldots,i_d),\,(j_1,\ldots,j_d)
\big)$.
Here, the mode index pair $(i_k,j_k)$ is formed by the
``(input, output)'' mode indices such that the application of $\ten{A}$
to $\ten{X}$, i.e., $\ten{A}\ten{X}$, transforms any input mode index
$j_k$ in $\ten{X}\big(\ldots,j_k,\ldots\big)$ into the output mode
index $i_k$ of $(\ten{A}\ten{X})\big(\ldots,i_k,\ldots\big)$, according to 
the formula
\begin{align*}
  \big(\ten{A}\ten{X}\big)(i_1,i_2,\dots,i_d)
  = \sum_{j_1=1}^{n_1}\sum_{j_2=1}^{n_2}\dots\sum_{j_d=1}^{n_d} \ten{A}(i_1,i_2,\dots,i_d,j_1,\dots,j_d)\ten{X}(j_1,\dots,j_d).
\end{align*}
To pursue this strategy further, we permute the mode indices of
$\ten{A}$ to pair together the input/output indices as in
\begin{equation}
  \ten{A}(i_1,i_2,\dots,i_d,j_1,j_2,\dots,j_d)
  \xrightarrow[\text{}]{\text{permute}}\ten{A}\big((i_1,j_1),(i_2,j_2),\dots,(i_d,j_d)\big).
  \label{eqn:permute}
\end{equation}
(We added the inner parenthesis to outline the index pairs).
The benefit of such permutation is that it helps separate/decompose
the dimensions whenever a ``matrix'' operator act on a single
dimension independently of the other dimensions.
The component-wise TT-format expression of the index-permuted, linear
operator $\ten{A}^{TT}$ reads as:
\begin{multline}
  \ten{A}^{TT}(i_1,j_1,i_2,j_2,\ldots,i_d,j_d)
  =  \sum^{r_1,r_2,\ldots,r_{d-1}}_{\alpha_{1},\alpha_2,\ldots,\alpha_{d-1}=1}
  \ten{A}_1\big(1,(i_1,j_1),\alpha_1\big)
  \ten{A}_2\big(\alpha_1,(i_2,j_2),\alpha_2\big) \ldots\\[0.5em]\ldots
  \ten{A}_d\big(\alpha_{d-1},(i_d,j_d),1\big),
  \label{eqn:TT-matrix-componentwise}
\end{multline}
where we make use of the 4D real cores
$\ten{A}_{k}=\big(\ten{A}_{k}(\alpha_{k-1},(i_k,j_k),\alpha_k\big)\in\mathbb{R}^{r_{k-1}\times m_k\times n_k\times r_{k}}$,
with $r_0=r_d=1$, and $k=1,2,\ldots,d$.
In~\eqref{eqn:TT-matrix-componentwise}, we enclosed the space indices
in parenthesis, i.e., $(i_k,j_k)$, to outline them.
Analogously to~\eqref{eqn:TT_def_element} and~\eqref{eqn:TT-vector},
we can reformulate \eqref{eqn:TT-matrix-componentwise} as the multiple
matrix product
\begin{equation}
  \ten{A}^{TT}(i_1,j_1,i_2,j_2,\ldots,i_d,j_d) =
  \mat{A}_1(i_1,j_1)
  \mat{A}_2(i_2,j_2) \dots
  \mat{A}_d(i_d,j_d).
  \label{eqn:TT-matrix-compact}
\end{equation}
\begin{figure}[h]
  \centering
  \includegraphics[width=\textwidth]{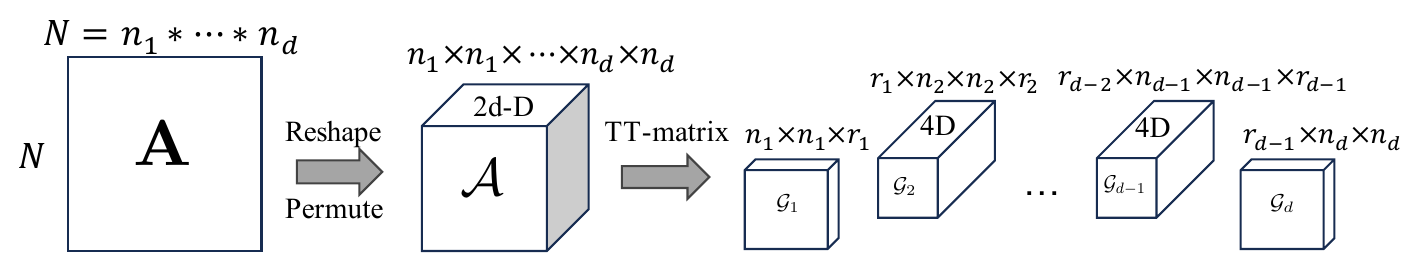}
  \caption{Representation of a linear operator in the TT-matrix
    format.
    First, we reshape the operation matrix $\mat{A}$ into a $2d$
    dimensional tensor $\ten{A}$ and permute its indices.
    Then, we factorize the tensor in the tensor-train matrix format
    according to Eq.~\eqref{eqn:TT-matrix-componentwise}.}
  \label{fig:TT-matrix-format}
\end{figure}
Equation~\eqref{eqn:TT-matrix-compact} is equivalent
to~\eqref{eqn:TT-matrix-componentwise} since each 4D core
$\ten{A}_{k}\big(\alpha_{k-1},(i_k,j_k),\alpha_{k}\big)$ can be seen
as an $r_{k-1}\times r_{k}$-sized, real matrix
$\big(\mat{A}(i_k,j_k)\big)_{\alpha_{k-1},\alpha_{k}}$ in the indices 
$\alpha_{k-1}$ and $\alpha_{k}$, where the index pair $(i_k,j_k)$ acts parametrically
for $i_k=1,2,\ldots,m_k$, $j_k=1,2,\ldots,n_k$.
Fig.~\ref{fig:TT-matrix-format} illustrates the steps needed to
construct the TT-matrix representation of a linear operator.

\medskip
We can further simplify the TT-matrix representations of $\ten{A}$ in
\eqref{eqn:TT-matrix-componentwise} and \eqref{eqn:TT-matrix-compact}
and the action of tensor $\ten{A}$ on tensor $\ten{X}$ in \TT{formats}
if the internal ranks of $\ten{A}$ are all equal to $1$.
In such a case, all summations in
Eq.~\eqref{eqn:TT-matrix-componentwise} reduce to a sequence of single
matrix-matrix multiplications, and $\ten{A}$ becomes the tensor
product of $d$ matrices:
\begin{align}
  \ten{A}^{TT} &= \mat{A}_1 \tenprod \mat{A}_2 \tenprod \dots \tenprod \mat{A}_d.
  \label{eqn:TT_matrices}
\end{align}
We show below how this approach works through the discretization of the 3D Laplace differential operator on a cubic domain. In the 1D case, a finite difference discretization can be represented
as a ${n\times n}$-sized, banded, Toeplitz matrix $\mat{L_1}$.
Assuming, for simplicity, that $n=n_1=n_2=n_3$, the $\big(n^{3}\times
n^{3}\big)$-sized matrix $\mat{L_3}$ of the Laplace operator can be
constructed as follows (see, e.g.,~\cite[Section
  3.1]{oseledets2011tensor}):
\begin{align*}
  \mat{L_3} =
  \mat{L_1} \otimes \mat{I_n} \otimes\mat{I_n} + 
  \mat{I_n} \otimes \mat{L_1} \otimes\mat{I_n} +
  \mat{I_n} \otimes\mat{I_n}  \otimes\mat{L_1},
\end{align*}
where $\mat{I_n}$ is the $n\times n$ identity matrix.
It is clear that $\mat{L}_3$ can be considered as a reshaping of the
$6D$ tensor $\ten{L}_3$ representing the Laplace differential operator
in the multi-dimensional setting.
Following the approach outlined above, we can directly construct
tensor $\ten{L}_3$ by replacing the Kronecker product $\otimes$ by the
tensor product $\tenprod$, to obtain that
\begin{align*}
  \ten{L}_3 =
  \mat{L}_1 \tenprod \mat{I_n} \tenprod \mat{I_n} +
  \mat{I_n} \tenprod \mat{L_1} \tenprod \mat{I_n} +
  \mat{I_n} \tenprod \mat{I_n} \tenprod \mat{L_1},
\end{align*}
and we can represent each of these three rank-1 terms in \TT{format}
using Eq.~\eqref{eqn:TT-vector}.
For example, the first term of the right-hand side has elements:
\begin{align*}
  \big(\mat{L}_1 \tenprod \mat{I_n} \tenprod \mat{I_n})^{TT}\big(i_1,j_1,i_2,j_2,i_3,j_3)
  =
  \mat{L}_1\big(1,i_1,j_1,1\big)
  \mat{I_n}\big(1,i_2,j_2,1\big)
  \mat{I_n}(1,i_3,j_3,1).
\end{align*}
Finally, we note that this process yields tensor $\ten{L}_3$ in
\TT{format} since the sum of tensors in \TT{format} is itself a tensor
in \TT{format}, see, e.g., \cite[]{oseledets2011tensor}, although a
rank reduction step could be necessary~\cite[]{oseledets2011tensor}.
Importantly, as we will show later, the discrete differential
operators of the NTEs possess similar
structures as the high-dimensional Laplace operator.

\begin{figure}[t]
  \centering
  \begin{tabular}{c}
    \includegraphics[width=0.75\textwidth]{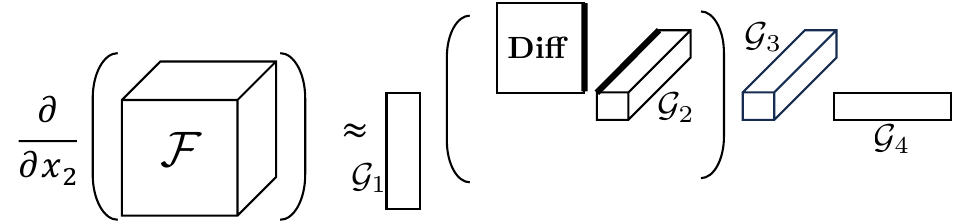}\\[1em]
    \includegraphics[width=0.75\textwidth]{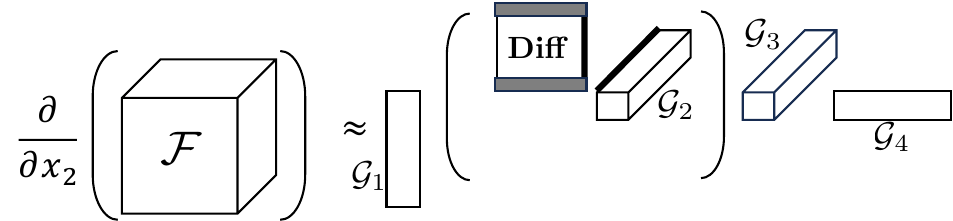}
  \end{tabular}
  \caption{
    Top panel: Computing derivative with respect to $x_2$ in \TT{format}
    by only applying the operation matrix $\mat{Diff}$ on the mode
    (\textbf{bold}) dimension of the second core $\ten{G}_2$.
    Bottom panel: Boundary conditions (grey) for the second dimension $x_2$
    are incorporated into the matrix operation of differentiation,
    $\mat{Diff}$.}
  \label{fig:TT_4D_x2_diff}
  \label{fig:TT_diff_x2_BC}
\end{figure}

\subsubsection{Boundary Conditions}
In general, incorporating boundary conditions into 
a discretization scheme involves modifying some entries 
of the discrete operators to account for the prescribed 
boundary conditions.
By separating the variables, the \TT{format} allows us to treat each
dimension independently, making it easier to handle boundary
conditions.
For example, the boundary conditions for the dimension $x_2$ can be
embedded into the differentiation matrix, e.g., $\mat{Diff}$ (see \ref{Appendix_B}, Eq.~\eqref{eq:diff:operator}), which, then,
is applied to the second core $\ten{G}_2$ as shown in
Fig~\ref{fig:TT_diff_x2_BC}.

\subsubsection{Quantized Tensor Train (QTT) Format}
\label{eqn:QTT_method}

\TT{format} is very effective for the compression of high-dimensional
tensors.
However, for large low-dimensional objects, such as vectors, matrices,
and tensors, the Quantized Tensor Train format \cite{khoromskij2011d}
is even more effective.
Specifically, we can reshape large vectors and matrices into
high-dimensional tensors with a small number of elements in each
dimension, and
subsequently decompose these tensors in the \TT{format}.

\subsubsection{Quantized Tensor Train representation of vectors}
\begin{figure}[httb]
  \centering
  \includegraphics[width=0.6\textwidth]{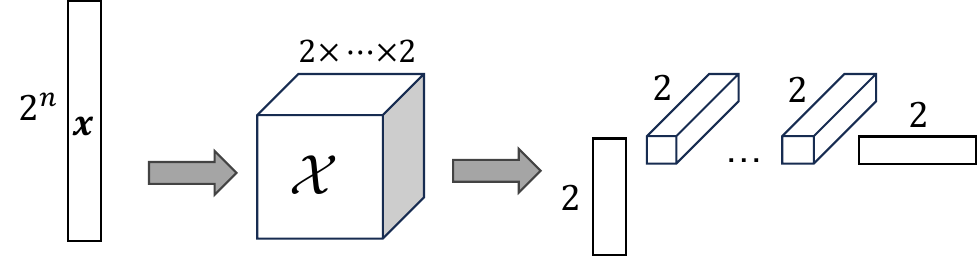}
  \caption{\QTT{format} of a vector $\mat{x}$. The vector is first
    reshaped into an $n$-dimensional tensor $\ten{X}$ of size
    $\bytwo$, and then the \TT{format} of that tensor is computed.}
    \label{fig:QTT_definition}
\end{figure}

Consider the vector $\mat{x} \in \mathbb{R}^{2^n}$.
We first reshape $\mat{x}$ to a $2\times2\times\ldots\times2$-sized
$n-$dimensional tensor, e.g.,
$\ten{X}\in\mathbb{R}^{2\times\ldots\times2}$.
Then, we decompose this tensor in \TT{format}, and called this
\TT{tensor representation} the \QTT{format}; see, e.g.,
Fig~\ref{fig:QTT_definition}.
The \QTT{format} exploits the fact that high-dimensional tensors, as
well as extra-large matrices that correspond to physics
\cite{udell2019big} often possess low-rank structures, meaning that the 
corresponding quantized tensor can be well approximated by a
combination of low-rank tensor cores.
QTT is an example of the so-called "blessing of
dimensionality" \cite{tyrtyshnikovblessing}, which is used to describe
the favorable properties that emerge when dealing with
high-dimensional data.
The ``blessing of dimensionality'' of QTT manifests itself in the form
of improved compression efficiency.

\subsubsection{Quantized Tensor Train representation of linear operators}
Let $\mat{A} \in \mathbb{R}^{2^n \times 2^n}$ be the matrix form of a
linear operator.
Converting this to a QTT matrix is a special case of the TT-matrix
format for linear operators, described in
Section~\ref{subsection:TT_operators}, where all the dimensions sizes
are 2.
First, we reshape matrix $\mat{A}$ into a $2n$-dimensional tensor
$\ten{A}$ of size $2\times2\times\ldots\times2$ ($n$ times).
Then, the dimensions of $\ten{A}$ is permuted as follows:
\begin{align*}
  \ten{A}(i_1,i_2,\ldots,i_n,j_1,\ldots,j_n) = \ten{A}\big( (i_1,j_1),(i_2,j_2)\ldots;(i_n,j_n)\big),
\end{align*}
finally, we decompose this tensor using the TT-matrix format
described in Eq.~\eqref{eqn:TT-matrix-componentwise} as
\begin{multline*}
  \ten{A}^{QTT}\big((i_1,j_1),(i_2,j_2),\ldots,(i_n,j_n)\big) =
  \sum^{r_1,\ldots,r_{d-1}}_{\alpha_{1},\ldots,\alpha_{d-1}=1}
  \ten{A}_1\big(1,i_1,j_1,\alpha_1\big)
  \ten{A}_2\big(\alpha_1,i_2,j_2,\alpha_2\big)
  \ldots\\
  \ldots
  \ten{A}_n\big(\alpha_{n-1},i_n,j_n,1\big),
\end{multline*}
where the first and the final cores are the 3D tensors
$\ten{A}_1\in\mathbb{R}^{2\times2\times r_1}$ and
$\ten{A}_{n}\in\mathbb{R}^{r_{l_k-1}\times 2\times 2}$, 
and each intermediate core is the 4D tensor
$\ten{A}_m\in\mathbb{R}^{r_{m-1}\times2\times 2 \times r_m}$, $m =
2\dots n-1$.


\section{Results}
\subsection{Tensorization of Boltzmann Neutron Transport Equation}
\label{sec4:tensorization-Boltzmann-NTE}



Numerically solving PDEs often leads to ultra-large linear algebra
problems, such as, linear systems of equations,
$\mat{A}\mat{x}=\mat{b}$, or generalized eigenvalue problems,
$\mat{A}\mat{x}=\lambda \mat{B}\mat{x}$.
To avoid the curse of dimensionality, we reformat all operators and
vectors in \TT{format} to achieve a computational complexity that grows
linearly, rather than exponentially,
\begin{equation*}
  \begin{aligned}
    \mat{A}\mat{x} &= \mat{b} &\rightarrow \ten{A}^{TT}\ten{X}^{TT} &= \ten{B}^{TT},\\
    \mat{A}\mat{x} &= \lambda \mat{C}\mat{x} &\rightarrow \ten{A}^{TT}\ten{X}^{TT} &=\lambda\ten{C}^{TT}\ten{X}^{TT}.
  \end{aligned}
\end{equation*}
Notice that $\ten{A}^{TT}$ and $\ten{C}^{TT}$ are in TT-matrix format (see subsection \ref{SUB:Lin_operinTT}),
while $\ten{X}$ and $\ten{B}$ are in \TT{format}.

If some of the \TT{cores} of $\ten{A}^{TT}$ or $\ten{C}^{TT}$ possess
special tensor structures, they can be further compressed.
For example, if a \TT{core} has a Toeplitz structure, we can compress
it using \QTT{format}, since Toeplitz matrices have low-rank
\QTT{formats} \cite{kazeev2012low}.
This is the case for the neutron transport problem we solve below.
Therefore, to achieve higher compression, we further transform these
large \TT{cores} into \QTT{format}, and then solve the mixed TT/QTT
version of the problems with an appropriate \TT{optimization}
technique (see \ref{APP:OPTTech}).
To find the biggest eigenvalue of the NTE problem we redesign a
fixed-point scheme to work with the \QTT{format}.

\medskip
Below, we outline the main steps of our method:

\begin{itemize}
\item \textbf{Forming the full tensor equation:} We form the
  discretization (that leads to linear system or generalized
  eigenvalue problem) with boundary conditions included in the
  discrete matricised form of the operators.
  To be able to apply the \QTT{format}, we choose the number of nodes of the grid
  in each dimension to be a power of two.

\item \textbf{\TT{format}:} We transform all objects in the full
  tensor equation into \TT{format}.
  
\item \textbf{\QTT{format}:} We transform the \TT{cores} that possess
  low-rank QTT structures into \QTT{format} to achieve higher
  compression.

\item \textbf{Solving the problem in \QTT{format}:} We apply the
  available \TT{solvers}, or design a new one, to solve the tensorized
  equations and obtain the solutions in \TT{format}.
\end{itemize}

Next, we show how to apply this general approach to solve the NTE.

\subsection{BNTE in \TT{format}}
\label{eqn:TT/QTT_metod}

In this section, we solve the discretization schemes for k-effective
(see Eq.~\eqref{eq:OpFormkNTE}) and alpha-eigenvalue criticality
problems (see Eq.~\eqref{eq:OpFormAlphaNTE}) utilizing the TT/QTT
format for BNTE.
In the traditional discretization, the operators $\mat{H}$, $\mat{S}$,
$\mat{F}$, and $\mat{V}^{-1}$ are designed as \textit{matrices} that
operate on the \textit{vectorized} solution $\Psi$.
This choice has been made to leverage existing linear algebra solvers
specifically tailored for matrices.
However, since the operation matrices are required to be fully formed,
the problem size is limited.
When the required storage exceeds the available memory,
the matrix-free approach is used.
However, using a matrix-free implementation can only help reducing the
memory usage but not the computational cost.

In our approach, we seek the compact \TT{format} of the eigenvector
$\Psi$ and the operators acting on it.
The eigenvector $\Psi$ in its original form is a five-dimensional
tensor with dimensions $G\times L\times K\times J\times M$.
Correspondingly, each operator acting on $\Psi$ can be represented as
a ten-dimensional tensor with dimensions
$(G\times G)\times(L\times L)\times(K\times K)\times(J\times
J)\times(M\times M)$.
We denote these tensor operators as $\ten{H},\ \ten{S},\ \ten{F},\ $
and $\ten{V}^{-1}$, where $\ten{H} = \ten{H}_x + \ten{H}_y + \ten{H}_z
+ \ten{H}_{\sigma}$.
Eqs. \eqref{eqn:term_for_operator} below lists the operators that we
need to construct in \TT{format} and their terms, which are given in
Equations \eqref{eq:Discretized_kNTE} and
\eqref{eq:Discretized_alphaNTE}:

\begin{equation}
  \begin{aligned}
    \frac{\mu_{\ell}}{4 \Delta x_{i}} \Bigg [ \sum_{k'=k-1}^{k} \sum_{j'=j-1}^{j} \psi_{g,\ell,k',j',i} - \psi_{g,\ell,k',j',i-1} \Bigg ]     &\rightarrow \ten{H}_x\Psi,\\
    \frac{\eta_{\ell}}{4 \Delta y_{j}} \Bigg [ \sum_{k'=k-1}^{k} \sum_{i'=i-1}^{i} \psi_{g,\ell,k',j,i'} - \psi_{g,\ell,k',j-1,i'} \Bigg ]    &\rightarrow \ten{H}_y\Psi,\\
    \frac{\xi_{\ell}}{4 \Delta z_{k}} \Bigg [ \sum_{j'=j-1}^{j} \sum_{i'=i-1}^{i} \psi_{g,\ell,k,j',i'} - \psi_{g,\ell,k-1,j',i'} \Bigg ]     &\rightarrow
    \ten{H}_z\Psi,\\ \frac{\sigma_{g,k,j,i}}{8} \Bigg [ \sum_{k'=k-1}^{k} \sum_{j'=j-1}^{j} \sum_{i'=i-1}^{i} \psi_{g,\ell,k',j',i'} \Bigg ]  &\rightarrow \ten{H}_\sigma\Psi,\\
    \frac{1}{k_{\text{eff}}}\frac{1}{8} \sum_{g'=1}^{G} \chi_{g'g}\nu\sigma_{f,g',k,j,i} \sum_{\ell' = 1}^{L} w_{\ell'} \Bigg [ \sum_{k'=k-1}^{k} \sum_{j'=j-1}^{j} \sum_{i'=i-1}^{i} \psi_{g',\ell',k',j',i'} \Bigg ] &\rightarrow \ten{S}\Psi,\\
    \frac{1}{8} \sum_{g'=1}^{G} \sigma_{s,g,g',k,j,i} \sum_{\ell' = 1}^{L} w_{\ell'} \Bigg [ \sum_{k'=k-1}^{k} \sum_{j'=j-1}^{j} \sum_{i'=i-1}^{i} \psi_{g',\ell',k',j',i'} \Bigg ] &\rightarrow \ten{F}\Psi, \\
    \frac{1}{8}\frac{\alpha}{v_{g}} \Bigg [ \sum_{k'=k-1}^{k} \sum_{j'=j-1}^{j} \sum_{i'=i-1}^{i} \psi_{g,\ell,k',j',i'} \Bigg ] &\rightarrow \ten{V}^{-1}\Psi.\\
  \end{aligned}
  \label{eqn:term_for_operator}
\end{equation}

These operators possess specific algebraic structures and we can
construct their TT-format representation explicitly through the
equivalence described in Eq.~\eqref{eqn:TT_matrices}.
Each of these operators contains several univariate operators,
i.e., operators acting on a single variable only.
This fact allows us to formulate them as tensor products of matrices,
which is equivalent to a \TT{format} with \TT{rank} equal to one.
As an example, we describe below the algebraic structure of one of
them, $\ten{H}_x$.
By definition, operator $\ten{H}_x$ explicitly depends on the five
indices $i$, $j$, $k$, $g$, $\ell$, which respectively discretize the
dimensions $x$, $y$, $z$, the energy dimension, and the angular
dimension.
The dependence on $i$ is through a first-order differentiation
formula; the dependence on $j$ and $k$ is through an average; the
dependence on $\ell$ is through a scaling factor, and the operator is
independent of $g$.
Specifically, 

\begin{itemize}
\item \textbf{Energy dimension:}~(index \textit{$g$}) - The operator
  $\ten{H}_x$ does not make any change in the energy dimension.
  Therefore, we can represent the action of this part of $\ten{H}_x$
  as the identity matrix $\mat{I}_{G}$.

\item \textbf{Ordinate dimension:}~(index \textit{$l$}) - For each
  $l$, the operator $\ten{H}_x$ multiply $\mu_l$ into the eigenvector
  $\Psi$.
  Therefore, we can represent this part of $\ten{H}_x$ as a block
  diagonal matrix, $\mat{Q}_{\mu}$, that contains all values of
  $\mu_l$.

\item \textbf{Interpolation along $z$ and $y$ dimensions:}~(indices
  \textit{$k$} and $j$) - Along these dimensions, the operator
  performs an average of $\Psi$ between any two adjacent nodes.
  This average is accomplished by the sums
  $\dfrac{1}{4}\sum_{k'=k-1}^{k} \sum_{j'=j-1}^{j}$.
  Therefore, we represent these parts of $\ten{H}_x$ by two
  interpolation matrices, $\intpmat_x$ and $\intpmat_y$.
  
\item \textbf{Differentiation along $x$ dimension:}~(index
  \textit{$i$}) - This part of the operator $\ten{H}_x$ performs a
  differentiation along $x$.
  Therefore, we represent this part of $\ten{H}_x$ through the
  differentiation matrix, $\mat{Diff}$.
  The boundary condition for the $x$ dimension dependents on whether
  $\mu_l$ is positive or negative, which requires a specific
  incorporation in the differentiation matrix,
  cf. Figure~\ref{fig:BC_1D}.
\end{itemize}

The remaining operators may depend differently on these indices, but the
way we treat them is similar.
In all cases we base their TT/QTT discretization on the tensor
operations discussed in Section~\ref{sec3:Tensor-Networks}.

This is the key observation that allow us to directly construct the
BNTE operators in \TT{format} as tensor products of univariate operators.

\subsubsection{\TT{format}s of all BNTE Left-Hand-Side (LHS) Operators}
The BNTE interaction tensor $\ten{H}$ that is on the LHS of the BNTE
includes four operation tensors, cf.~Eq.~\eqref{eq:Discretized_kNTE}:
\begin{align*}
  \ten{H} = \ten{H}_x + \ten{H}_y + \ten{H}_z + \ten{H}_{\sigma}.
\end{align*}
Below we explicitly construct the operation matrices that are the
\TT{cores} in the \TT{format} of the interaction operator $\ten{H}$.

\begin{wrapfigure}{R}{0.3\textwidth}
  \centering
  \includegraphics[width = 0.3\textwidth]{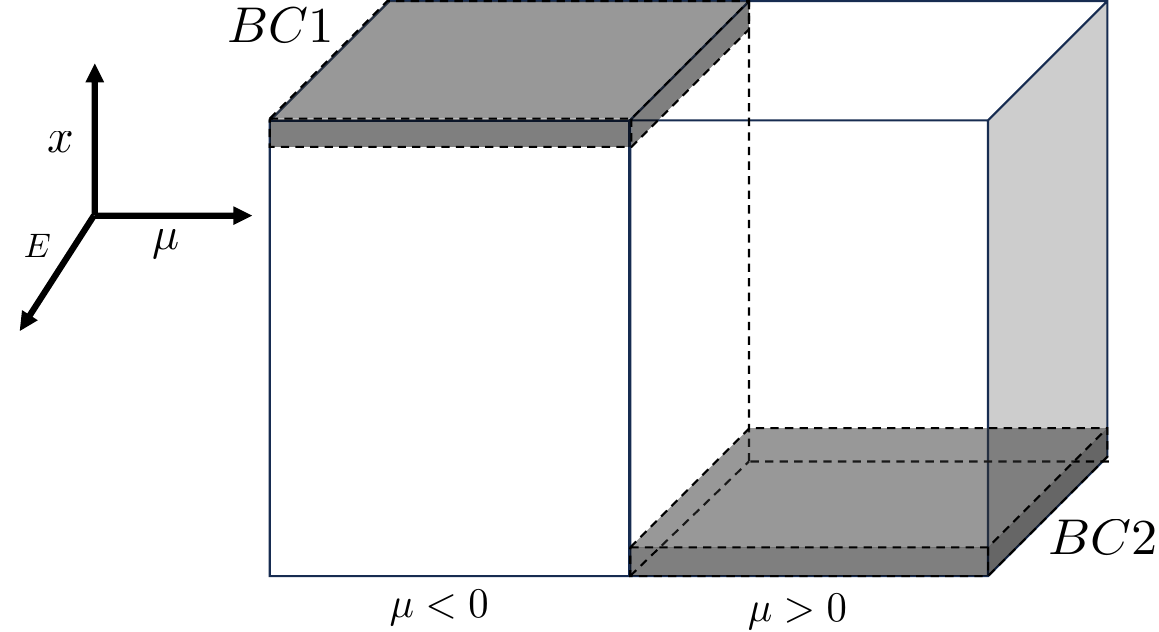}
  \caption{Boundary conditions BC1 and BC2 dependent on either $\mu$'s
    values are positive or negative.}
  \vspace{-2.cm}
  \label{fig:BC_1D}
\end{wrapfigure}

\paragraph{Differentiation Matrices}
The differentiation matrix, $\mat{D}_x$, is used to compute the finite
difference along spatial dimensions $x$.
The boundary conditions depend on the sign of the ordinate, see
Eq.~\eqref{eq:discreteBCs}.
In Fig. \ref{fig:BC_1D} we show an example of two boundary conditions
for the 1D neutron transport equation, which has three variables,
i.e. $x$, ordinate $\mu$, and energy $E$.
When $\mu$ is positive, the boundary condition BC1 is at $i=M$, when $\mu$ is negative, the boundary condition BC2 is at $i=0$.

For the positive ordinates $\mu_\ell$, the differentiation matrix
acting on $x$ dimension, $\diffmat^+_x$ is defined as:
\begin{equation}
  \diffmat^+_x \equiv
  \setlength\arraycolsep{2pt}
  \dfrac{1}{\Delta x}\begin{pmatrix}
    1  &        &        &\\
    -1 & \ddots &        & \\
    {} & \ddots & \ddots &\\
    {} &        & -1     & 1
  \end{pmatrix} \in \mathbb{R}^{(M+1) \times (M+1)},
\end{equation}
For the negative ordinates $\mu_\ell$, the differentiation matrix
acting on $x$ dimension, $\diffmat^-_x$ is defined as:
\begin{equation}
  \diffmat^-_x \equiv
  \setlength\arraycolsep{2pt}
  \dfrac{1}{\Delta x}\begin{pmatrix}
    -1 & 1 & &\\
    & \ddots & \ddots &\\
    & & \ddots & 1\\
    & & & -1
  \end{pmatrix} \in \mathbb{R}^{(M+1) \times (M+1)}.
\end{equation}
In this way, we explicitly incorporate the boundary conditions in the
differentiation matrices.
The construction of the differentiation matrix in 3D for $y$ and $z$,
i.e. $\diffmat^+_y$, $\diffmat^-_y$, $\diffmat^+_z$, and $\diffmat^-_z$,
is the same.

\paragraph{Interpolation Matrices}
We use the interpolation matrix to approximate the BNTE solution, i.e.,
$\mat{\Psi}$, in the centers of the grid cells based on vertex values.
Similar to the differentiation matrix, the interpolation matrices
incorporate the boundary conditions depending on the sign of the
ordinate values.
The interpolation matrices acting on the $x$ dimension for $\mu_\ell
<0$, $\intpmat_x^-$, and for $\mu_\ell >0$, $\intpmat_x^+$, are:
\begin{equation}
  \intpmat^-_x \equiv
  \setlength\arraycolsep{3pt}
  \dfrac{1}{2}\begin{pmatrix}
    1  &      1 &        &  \\
    {} & \ddots & \ddots &  \\
    {} &        & \ddots & 1\\
    {} &        &        & 1
  \end{pmatrix} \in \mathbb{R}^{(M+1) \times (M+1)},
\end{equation}

\begin{equation}
  \intpmat^+_x \equiv
  \setlength\arraycolsep{3pt}
  \dfrac{1}{2}\begin{pmatrix}
    1  &        &        &\\
    1  & \ddots &        &\\
    {} & \ddots & \ddots &\\
    {} &        &      1 & 1
  \end{pmatrix} \in \mathbb{R}^{(M+1) \times (M+1)}.
\end{equation}

We also use interpolation matrices that do not include the boundaries:

\begin{equation}
  \intpmat^-_{x,noBC} \equiv
  \setlength\arraycolsep{3pt}
  \dfrac{1}{2}\begin{pmatrix}
    1  & 1   & &\\
    {} & \ddots & \ddots &\\
    {} &     & 1 & 1\\
    {} &     & & 0
  \end{pmatrix} \in \mathbb{R}^{(M+1) \times (M+1)}.
\end{equation}

\begin{equation}
  \intpmat^+_{x,noBC} \equiv
  \setlength\arraycolsep{3pt}
  \dfrac{1}{2}\begin{pmatrix}
    0  &        &        &\\
    1  &      1 &        &\\
    {} & \ddots & \ddots &\\
    {} &        &      1 & 1
  \end{pmatrix} \in \mathbb{R}^{(M+1) \times (M+1)},
\end{equation}

\paragraph{Angular Point Matrices}
The angular matrices are used to multiply appropriate angular points
to the correct terms on the left-hand-side of
Eqs.~\eqref{eq:Discretized_kNTE} and \eqref{eq:Discretized_alphaNTE}.
The dimensions of the angular matrices are $L \times L$.
The angular matrices are 8-block diagonal matrices, in which each
block corresponds to one of eight boundary conditions
in~\eqref{eq:discreteBCs}.
For each boundary condition, there are three angular matrices
\{$\angmat_{\mu},\angmat_{\eta},\angmat_{\xi}$\}, one for each angle
$\mu$, $\eta$, and $\xi$. For example, given the first boundary
condition $\mu_\ell<0, \eta_\ell<0, \xi_\ell<0$, the angular matrices
are formed as follows:

\begin{equation}
  \begin{cases}
    \mathbf{Q}_{\mu} = \mathbf{C}_{1} \otimes \text{diag}(\hat{\mu}_{-})\\[0.5em]
    \mathbf{Q}_{\eta} = \mathbf{C}_{1} \otimes \text{diag}(\hat{\eta}_{-})\\[0.5em]
    \mathbf{Q}_{\xi} = \mathbf{C}_{1} \otimes \text{diag}(\hat{\xi}_{-})
  \end{cases},\quad \mathbf{C}_{1} \equiv
  \setlength\arraycolsep{3pt}
  \begin{pmatrix}
    1  &   &        &\\
    {} & 0 &        &\\
    {} &   & \ddots &\\
    {} &   &        & 0
  \end{pmatrix} \in \mathbb{R}^{8 \times 8},
\end{equation}

where $\text{diag}(\mu)$ is a square diagonal matrix with the elements
of vector $\mu$ on the main diagonal; $\mathbf{C}_{i}$ is a $8\times8$
matrix with $\left(\mat{C}_{i}\right)_{i,i} = 1$, and zeros elsewhere;
index $i$ reflects the boundary condition.
The quantities $\hat{\mu}_{-}, \hat{\eta}_{-}$, and $\hat{\xi}_{-} $
are defined in
Section~\ref{subsubsec:discrete-ordinates-approximation}.

%

\paragraph{Boundary Conditions}
The $\ten{H}_j^{TT,bci}$ for the boundary condition $bci$, and for the
variable $j \in \{x,y,z,\sigma\}$ can be explicitly built using the
operation matrices defined above.
For example, the $\ten{H}_j^{TT,bci=1}$, for the boundary condition 1
( $\mu_\ell<0, \eta_\ell<0, \xi_\ell<0$ ) is defined as follows:\\
\begin{equation}
  \begin{aligned}
    &\ten{H}_x^{TT,1}=  \mathbf{I}_{G}  \tenprod \angmat_\mu \tenprod \intpmat_z^- \tenprod \intpmat_y^- \tenprod  \diffmat_x^-, \\[0.5em]
    &\ten{H}_y^{TT,1}=  \mathbf{I}_{G}\tenprod \angmat_\eta \tenprod   \intpmat_z^- \tenprod \diffmat_y^-  \tenprod \intpmat_x^-,\\[0.5em]
    &\ten{H}_z^{TT,1}= \mathbf{I}_{G} \tenprod \angmat_\xi \tenprod  \diffmat_z^- \tenprod \intpmat_y^- \tenprod \intpmat_x^-,\\[0.5em]
    &\ten{H}_\sigma^{TT,1}= \text{diag}(\sigma) \tenprod (\mathbf{C}_{1}\otimes\mathbf{I}_{L/8})\tenprod \intpmat_z^- \tenprod \intpmat_y^-  \tenprod \intpmat_x^-. 
  \end{aligned}
\end{equation}
Other  $\ten{H}_j^{TT,bci}$ are similarly constructed.

\paragraph{Velocity Tensor $\ten{V}^{-1}$}
For solving the alpha-eigenvalue problem, we need the velocity tensor
operator $\ten{V}^{-1}$, cf. Eq.~\eqref{eq:Discretized_alphaNTE}.
Similar to other operators, we construct the TT-format representation of $\ten{V}^{-1}$ from its \TT{format} for each boundary condition.
For the boundary condition labeled by \textit{bci}, the \TT{format} representation of $\left(\ten{V}^{-1}\right)^{TT, bci}$ is:\\
\begin{equation}  
  \left(\ten{V}^{-1}\right)^{TT, bci} = \text{diag}(1/\mat{v})
  \tenprod (\mathbf{C}_{1}\otimes\mathbf{I}_{L/8}) \tenprod \intpmat_z^- \tenprod \intpmat_y^- \tenprod  \intpmat_x^-.
\end{equation}
where $\mat{v}$
is the velocity vector of all energy groups, and $1/\mat{v}$ is element-wise division.

\subsubsection{\TT{format}s of the BNTE Righ-Hand-Side (RHS) Operators}
In view of Eq.~\eqref{eq:Discretized_kNTE}, the RHS of the BNTE 
consists of two tensor operators,  $\ten{S}$ and $\ten{F}$.
In the next two subsections, we explicitly construct the operation matrices 
that we use  as \TT{cores} in the \TT{format} representation of these 
operators.

\paragraph{Integral Operator Matrix}
Given $w \in \mathbb{R}^{L/8}$ is the weight vector
from~\eqref{eq:quadrature points}.
The integral operator matrix for the first boundary condition is:
\begin{align*}
  \mathbf{Intg}^1 = \mathbf{C}_{1}\otimes (\mathbf{1}_{L/8} \otimes w )
\end{align*}
where $\mathbf{1}$ is a vector of all ones.\\

\paragraph{Constructing $\ten{F}^{TT}$ and $\ten{S}^{TT}$}
Each operator $\ten{F}^{TT,bci}$ can be constructed using the above
defined operation matrices.
The $\ten{F}^{TT,1}$ matrix, for the first boundary condition, is
shown below as an example:
\begin{equation}
  \begin{aligned}
    &\ten{F}^{TT,1} =  \text{diag}(\nu\sigma_f) \tenprod\mathbf{Intg}^1 \tenprod \intpmat_{z,noBC}^- \tenprod  \intpmat_{y,noBC}^- \tenprod  \intpmat^-_{x,noBC}.\\
  \end{aligned}
\end{equation}
The operator $\ten{S}^{TT,1}$ can be similarly constructed as follows:
\begin{equation}
  \begin{aligned}
    &\ten{S}^{TT,1} = \text{diag}(\sigma_s) \tenprod\mathbf{Intg}^1 \tenprod \intpmat_{z,noBC}^- \tenprod  \intpmat_{y,noBC}^- \tenprod  \intpmat^-_{x,noBC}\\
  \end{aligned}
\end{equation}
For the derivation of the remaining $\ten{F}^{TT,bci}$,
$\ten{S}^{TT,bci}$, we refer to~\ref{eqn:Appendix}.
Then, $\ten{F}^{TT,bci}$ and $\ten{S}^{TT,bci}$ are used to construct
$\ten{F}^{TT}$ and $\ten{S}^{TT}$ :
\begin{equation}
  \begin{aligned}
    \ten{F}^{TT} = \sum_{bci=1}^8 \ten{F}^{TT,bci},~ 
    \ten{S}^{TT} = \sum_{bci=1}^8 \ten{S}^{TT,bci}.
  \end{aligned}
\end{equation}

\subsection{Transforming the BNTE Operators from TT to \QTT{format}}

Now that we have BNTE operators in \TT{format} with \TT{rank} one, we
present a general strategy to reformat them in QTT format, using the
fact that the operation matrices of the \TT{format}ted tensors have a 
Toeplitz structure, \cite{kazeev2013multilevel}.
To use the QTT format, the dimension sizes
$G,\ L, \ K,\ J,$ and $M $ must be a power of two. 
It can be easily
achieved by discretizing the full tensor by choosing the
number of nodes in each dimension to be a power of two.


Given that all operators
$\ten{H}_x,\ \ten{H}_y,\ \ten{H}_z,\ \ten{H}_{\sigma}, \ten{S},
\ten{F}$, and $\ten{V}^{-1}$ possess similar TT structures, we
describe the procedure to construct the \QTT{format} for
the generic operator $\ten{A} \in
\big\{\ten{H}_x,\ \ten{H}_y,\ \ten{H}_z,\ \ten{H}_{\sigma}, \\ \ten{S},
\ten{F}, \ten{V}^{-1} \big\}$ .

For every boundary condition, indexed by $bci \in \{1,\dots,8\}$ (see
Appendix and Eq.~\eqref{eq:discreteBCs}), there exists a TT-matrix format
representation of the generic discrete operator $\ten{A}^{TT}$,
with \TT{rank} one, which is specific to that
boundary condition, hereafter denoted by $\ten{A}^{TT,bci}$,
\begin{align*}
  \ten{A}^{TT,bci} = \mat{G}_1\tenprod \mat{G}_2\tenprod \mat{G}_3 \tenprod \mat{G}_4 \tenprod \mat{G}_5.
\end{align*}
The next step is to convert the above operation matrices $\mat{G}_k$
into their corresponding QTT-matrix formats.
These matrices are of Toeplitz structures. The procedure to convert
$\ten{A}^{TT,bci}$ to $\ten{A}^{QTT,bci}$ is described in the
Algorithm \ref{alg:TTtoQTT}.

\begin{algorithm}[H]
  \caption{Convert $\ten{A}^{TT}$ into $\ten{A}^{QTT}$}
  \label{alg:TTtoQTT}
  \KwData{$\ten{A}^{TT}=\mat{G}_1\tenprod \mat{G}_2 \tenprod \mat{G}_3 \tenprod \mat{G}_4 \tenprod \mat{G}_5$\\
    Matrix $\mat{G}_k$ has dimension of $n_k \times n_k$, where $n_k = 2^{l_k}$}
  \KwResult{$\ten{A}^{QTT}$}
  
  \For{$k = 1:5$}{
    Reshape $\mat{G}_k$ to $\bytwo$ tensor $\ten{G}_k$\\
    Permute $\ten{G}_k(i_1,\ldots,i_{l_k},j_1,\ldots,j_{l_k})$ to $\ten{G}_k(i_1,j_1,\ldots,i_{l_k},j_{l_k})$\\
    Compute $\ten{G}_k$'s QTT-matrix format $\ten{G}^{QTT}_k = g_1g_2 \dots g_{l_k}$\\
    \nonl \ where the first core is 3D tensor $g_1 \in
    \mathbb{R}^{2\times2\times r_1}$, the final core is a 3D tensor
    $g_{l_k} \in \mathbb{R}^{r_{l_k-1}\times 2\times 2}$ and each
    middle core, $g_m$, is a 4D tensor, $g_m \in
    \mathbb{R}^{r_{m-1}\times2\times 2 \times r_m}$ for $m = 2 \dots
    l_k-1$
  }
  $\ten{A}^{QTT} = \ten{G}^{QTT}_1\ten{G}^{QTT}_2\ten{G}^{QTT}_3\ten{G}^{QTT}_4\ten{G}^{QTT}_5$
\end{algorithm}

Finally, the \QTT{format} of the operator, $\ten{A}^{QTT}$, is the sum
of all \QTT{format} of that operator for each boundary condition:
\begin{align*}
  \ten{A}^{QTT} = \sum_{bci=1}^8 \ten{A}^{QTT,bci}.
\end{align*}

\subsubsection{\QTT{format} of the Interaction Tensor $\ten{H}^{QTT}$}
In Sec. \ref{eqn:QTT_method} we described the general procedure to
construct the \QTT{format} for an operator.
Importantly, we show above that the BNTE operators in \TT{format} are
tensor products of univariate operation matrices with specific
algebraic structure.
Therefore, next we define these matrices (or \TT{cores}, $\mat{G}_k$) to
construct each one of the NTE operators in the list
$\{\ten{H}_x,\ \ten{H}_y,\ \ten{H}_z,\ \ten{H}_{\sigma}, \ten{S},
\ten{F}, \ten{V}^{-1} \}$ and transform them in \QTT{format} exploiting
their algebraic structure.
We will start by describing how to construct the \QTT{format}.
Now we have all needed matrices to obtain the interaction operator
$\ten{H}$ in \QTT{format}, which we will denote by, $\ten{H}^{QTT}$.
Then using Algorithm \ref{alg:TTtoQTT}, $\ten{H}_j^{TT,bci}$ is
converted into the \QTT{format} as $\ten{H}_j^{QTT,bci}$.
Finally, the \QTT{format} of the interaction tensor $\ten{H}$ can be
computed as:
\begin{equation}
  \ten{H}^{QTT} = \sum_{j \in \{x,y,z,\sigma\}} \sum_{bci=1}^8 \ten{H}_j^{QTT,bci}
\end{equation}


\subsubsection{\QTT{Format} of Fission Operator $\ten{F}^{QTT}$ and Scattering Operator $\ten{S}^{QTT}$}
The QTT representations of fission operator $\ten{F}^{QTT}$ and
scattering operator $\ten{S}^{QTT}$ are constructed in a similar
fashion as $\ten{H}^{QTT}$.
For each boundary condition, the component matrices are identified and
converted to QTT format before being merged to finalized the QTT
representation of the operator for that boundary condition.

\subsubsection{Fixed-Point Algorithms in \QTT{format}}
We have completed constructing the \QTT{format} of all operators needed
to rewrite Eqn.~\eqref{eq:Discretized_kNTE} and
\eqref{eq:Discretized_alphaNTE} into \QTT{format}.
Next we will describe the \TT{solvers}, and the fixed-point method \cite{shashkin1991fixed} we utilize here, 
we use to solve these tensorized equations.

\subsubsection{K-Effective Problem}
A \TT{format} analogous to Eq.~\eqref{eqn:k_eff_fixed_points}
approximates the eigenvalue $k_\text{eff}$ and the eigenvector
$\Psi^{QTT}$ in the \QTT{format} for the following problem:
\begin{equation}
  \ten{H}^{QTT} \Psi^{QTT} = \left[ \ten{S}^{QTT} + \dfrac{1}{k_\text{eff}\ten{F}^{QTT}}\right]\Psi^{QTT} \
\end{equation}
Starting from random initial guesses for $\Psi^{QTT}$ and
$k{_\text{eff}}$, this equation can be solved using the following
fixed point scheme:
\begin{equation}
  \begin{aligned}
    \ten{B}^{QTT,(\tau)}& =\bigg [ \ten{S}^{QTT} + \frac{1}{k_{\text{eff}}^{(\tau)}} \ten{F}^{QTT} \bigg ] \Psi^{QTT,(\tau)}\\
    \ten{H}^{QTT}\Psi^{QTT,(\tau+1)} &= \ten{B}^{QTT,(\tau)}\ \\
    k_{\text{eff}}^{(\tau+1)} &= k_{\text{eff}}^{(\tau)} \frac{ \sum \ten{F}^{QTT} \Psi^{QTT,(\tau+1)}}{ \sum \ten{F}^{QTT} \Psi^{QTT,(\tau)}}, 
  \end{aligned}
  \label{eqn:TT_k_eff_fixed_points}
\end{equation}
where the last update rule comes from the traditional matrix-free
approach, Eq.~\eqref{eqn:TT_k_eff_fixed_points}.

\subsubsection{Alpha-eigenvalue Problem}
Similarly, a \TT{format} analog of
Eq.~\eqref{eqn:alpha_eig_fixed_points} approximates the eigenvalue
$\alpha$ and the eigenvector $\Psi^{QTT}$ for the problem:
\begin{equation}
  \left[ \alpha\ten{V}^{-1,QTT} +  \ten{H}^{QTT} \right] \Psi^{QTT} = \left[ \ten{S}^{QTT} +\ten{F}^{QTT} \right]\Psi^{QTT} \
\end{equation}
Starting from random values for $\alpha^{(0)} = 0$, $\alpha^{(1)} =
0.01$, the alpha-eigenvalue problem in
Eq.~\eqref{eq:Discretized_alphaNTE} can be solved as follows.
First, we prepare the initial conditions by solving two
$k_{\text{eff}}$ problems, starting from random guesses for
$\Psi^{QTT}$ and $k_{\text{eff}}$ and using updating rules given by
Eq.~\eqref{eqn:TT_k_eff_fixed_points}:
\begin{equation}
  \begin{aligned}
    \bigg [ \alpha^{(0)}\ten{V}^{-1,QTT} + \ten{H}^{QTT} \bigg ]\Psi^{QTT,(0)} &= \bigg [ \ten{S}^{QTT} +  \frac{1}{k_{\text{eff}}^{(0)}} \ten{F}^{QTT} \bigg ] \Psi^{QTT,(0)}\\
    \bigg [ \alpha^{(1)}\ten{V}^{-1,QTT} + \ten{H}^{QTT} \bigg ]\Psi^{QTT,(1)} &= \bigg [ \ten{S}^{QTT} +  \frac{1}{k_{\text{eff}}^{(1)}} \ten{F}^{QTT} \bigg ] \Psi^{QTT,(1)},\\
  \end{aligned}
\end{equation}
which results in the needed four initial conditions: $\alpha^{(0)} =
0$, $\alpha^{(1)} = 0.01$, $k_{\text{eff}}^{(0)}$, and
$k_{\text{eff}}^{(1)}$.
Then using the update rule from traditional matrix-free approach,
Eq.~\eqref{eqn:alpha_eig_fixed_points}, we can calculate
$\alpha^{(j+1)}$,
\begin{equation}
  \begin{aligned}
    \alpha^{(\tau+1)} &= \alpha^j + \dfrac{1-{k_{\text{eff}}^{(\tau-1)}}}{(k_{\text{eff}}^{(\tau)}-k_{\text{eff}}^{(\tau-1)})(\alpha^{(\tau)} - \alpha^{(\tau-1)})}.\\
  \end{aligned}
\end{equation}
Second, we solve the next $k_{\text{eff}}$ problem with the initial
conditions, $\alpha^{(\tau+1)}$, and again using the updating rules given
by Eq.~\eqref{eqn:TT_k_eff_fixed_points} with random $k_{\text{eff}}$,
$\Psi^{QTT}$:
\begin{equation}
  \begin{aligned}
    \bigg[ \alpha^{(\tau+1)}\ten{V}^{-1,QTT} + \ten{H}^{QTT} \bigg]\Psi^{QTT,(\tau+1)} &= \bigg[\ten{S}^{QTT} +  \frac{1}{k_{\text{eff}}^{(\tau+1)}} \ten{F}^{QTT} \bigg ] \Psi^{QTT,(\tau+1)}
    \label{eqn:TT_alpha_eig_fixed_points}
  \end{aligned}
\end{equation}

\subsection{Tensor Optimization Techniques}
\label{APP:OPTTech}
In order to apply the update rules (described in the previous section)
we need to be able to perform matrix/vector multiplications and solve
systems of linear equations in QTT format.
These tasked can be performed very efficiently in \TT{format}, by
manipulating a single \TT{core} at the time, and using tensor
optimization techniques,
see~\cite{holtz2012alternating,savostyanov2011fast,dolgov2014alternating}.

These linear algebra problems first need to be formulated as
optimization problems.
For example, a linear system $\mat{A}\mat{x} = \mat{b}$ where
$\mat{x}$, $\mat{b} \in \mathbb{R}^{N}$ and $\mat{A} \in
\mathbb{R}^{N\times N}$, can be solved as the following minimization
problem:
\begin{equation}
  \mat{x} = \argmin_{\mat{y}} \mathcal{J}(\mat{y}) \text{ where }  \mathcal{J}(\mat{y})
  = \lVert \mat{A}\mat{y} - \mat{b} \rVert^2
  = \mat{y}^T\mat{A}^T\mat{A}\mat{y} - 2\mat{b}^T\mat{A}\mat{y} + \mat{b}^T\mat{b}.
  \label{eqn:linear_system_obj}
\end{equation}
Similarly, a matrix-vector multiplication $\mat{A}\mat{b}$ can be
formulated as:
\begin{equation}
  \mat{x} = \argmin_{\mat{y}} \mathcal{J}(\mat{y}) \text{ where }  \mathcal{J}(\mat{y})
  = \lVert \mat{A}\mat{b} - \mat{y} \rVert^2
  = \mat{b}^T\mat{A}^T\mat{A}\mat{b} - 2\mat{y}^T\mat{A}\mat{b} + \mat{y}^T\mat{y}.
  \label{eqn:mv_obj}
\end{equation}

When the size of $\mat{A}$ and $\mat{b}$ are small enough, performing
the matrix-vector multiplication $\mat{A}\mat{b}$ directly is
straightforward.
However, when their sizes are too big, and they possess low-rank TT
structures, \TT{format} is obviously a better option.
In the \TT{format}, there are an explicit and exact algorithm to compute
the matrix-vector multiplication \cite{oseledets2011tensor}, but when
compared to the tensor optimization techniques, they perform much
slower.
This is the reason we choose to perform the matrix-vector
multiplication in \TT{format} using the tensor optimization techniques.

For both problems in Eqs. \eqref{eqn:linear_system_obj} and
\eqref{eqn:mv_obj}, when $\mat{A}$, $\mat{b}$, and $\mat{y}$ are in
\TT{format} as $\ten{A}^{TT}$, $\ten{B}^{TT}$, and $\ten{Y}^{TT}$,
$\mat{x}$ is also in \TT{format}, and:
\begin{equation}
  \ten{X}^{TT} = \argmin_{\ten{Y}^{TT}} \mathcal{J}(\ten{Y}^{TT})
  \label{eqn:TT_min_LS_problem}
\end{equation}
Formulating the optimization problem in \TT{format} enables the usage of
algorithms that fix all but one \TT{core} of $\ten{Y}^{TT}$, and turn
the multilinear problem into a series of much smaller linear problems
for each \TT{core}.
These algorithms include \emph{Alternating Linear Scheme (ALS)}
\cite{holtz2012alternating}, \emph{Two-Site Density Matrix
Renormalization Group (DMRG)} \cite{savostyanov2011fast}, or
\emph{Alternating Minimal Energy (AMEn)} \cite{dolgov2014alternating},
that can solve the minimization problem and find the optimal \TT{rank}
for it, \textbf{without ever working with the full tensor}.
In general, these algorithms need as an input, $\ten{A}^{TT}$ and
$\ten{B}^{TT}$, as well as the initial guess for $\ten{X}^{TT}$
(usually a random tensor with some \TT{rank}), and the acceptable
error, $\varepsilon$, for the solution.
We have used \texttt{amen\_solve}, and \texttt{amen\_mv} in the MATLAB
TT-Toolbox \cite{tt-toolbox} to perform these calculations.


\subsection{Numerical Experiments}
\label{sec5:numerical-experiments}



\subsubsection{A One-Dimensional Case Study}

A one-dimensional slab problem was considered to verify the correctness of our TT-method and measure its performance. This problem is part of a criticality verification benchmark suite \cite{sood2003} and is exactly critical ($k_{\text{eff}}=1$, $\alpha = 0)$. It consists of a plutonium-239 slab. The problem cross section data is shown in Table \ref{tab:1dbench}. The problem has energy one-group and the slab width is 3.707444 $cm$.

\begin{table}[httb]
    \centering

\begin{tabular}{|c|c|c|c|c|}
\hline
Material &  $\nu$ & $\sigma_{f}$ $(cm^{-1})$ & $\sigma_{s}$ $(cm^{-1})$ & $\sigma_{t}$ $(cm^{-1})$ \\
\hline
Pu-239 & 3.24 & 0.081600 & 0.225216 & 0.32640 \\
\hline
\end{tabular}


    \caption{Cross Section Data for 1D Benchmark Problem}
    \label{tab:1dbench}
\end{table}

We compared three approaches to solve this k-effective problem. The first approach uses a standard generalized eigen-solver (\textbf{GES}) to solve for the $k_{\text{eff}}$. 
The second approach is the iterative solver described in equation \ref{eqn:k_eff_fixed_points}applied to the full matricized format of the operators (\textbf{ISFM}), 
and the third approach is our TT/QTT iterative solver we build (\textbf{ISTT}). 
For the parameters, we set the number of spacial point to be 1024, and vary the number of ordinates $L \in \{2,4,6,8,16,32\}$. These sizes are small enough to show how accurate the approximations are for both eigenvalues and eigenvectors. Moreover, given that we know the ground truth values for the $k_{\text{eff}}$, we can show that by increasing $L$, the approximated eigenvalue should converge to the ground truth, ($k_{\text{eff}}=1$).
For the iterative solvers and for the tensor train solver, the tolerance was set at $10^{-6}$.

\begin{figure}[httb]
    \centering
    \includegraphics[width=1.0\textwidth]{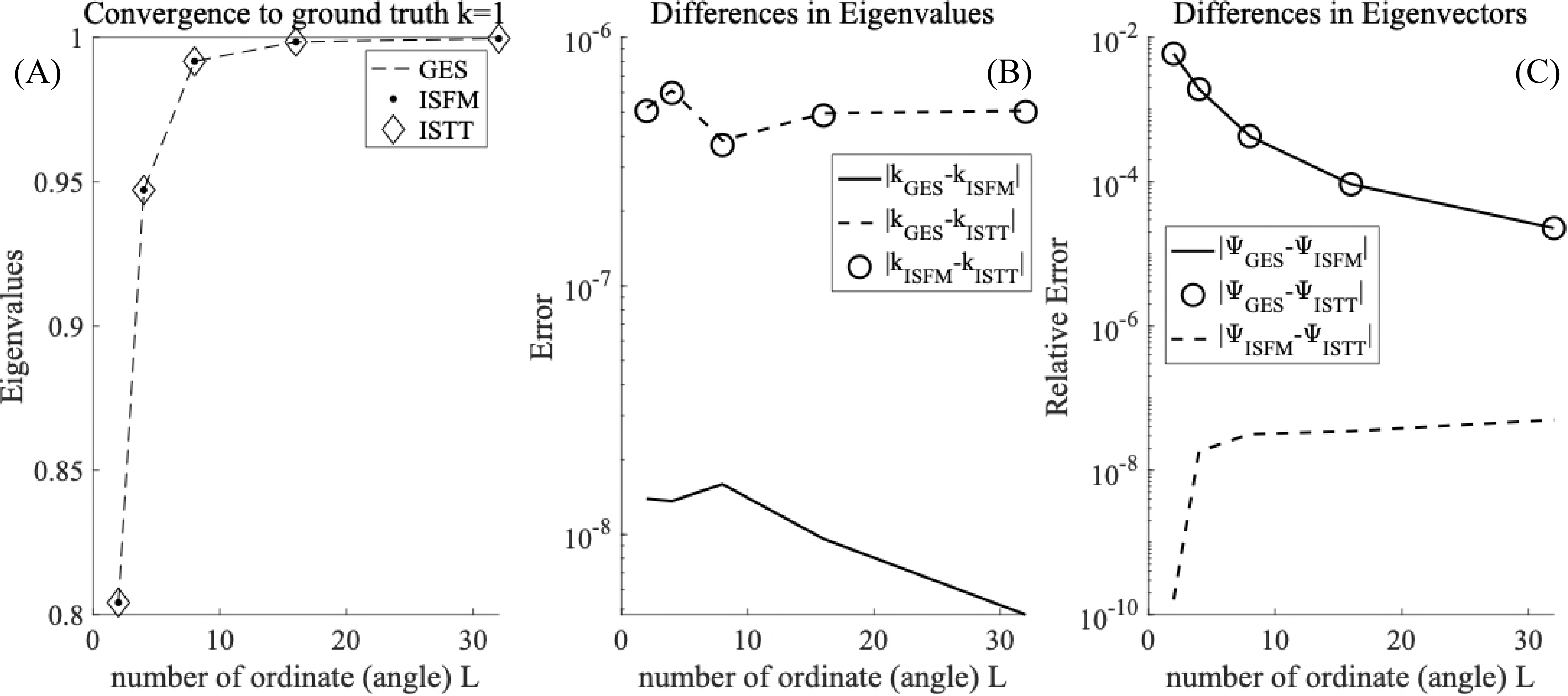}
    \caption{Comparison between three approaches to solve the k-effective problem. \textbf{GES} uses the generalized eigenvalue solver. \textbf{ISFM} is the iterative method for full matrix format of the operators. \textbf{ISTT} denotes our TT approach. (\textit{Panel A}) - Approximated eigenvalues converge to the ground truth value as the number of ordinate increases. (\textit{Panel B}) - Differences between approximated eigenvalues shows that the tensor train solver can approximate the eigenvalues with high accuracy. (\textit{Panel C}) - Differences between approximated eigenvectors shows that the eigenvector approximated by the tensor train solver is well matched with the ones from the full matrix iterative solver.}
    \label{fig:Pu239_1}
\end{figure}

Figure \ref{fig:Pu239_1} shows how well the eigenvalues and eigenvectors are approximated. Figure \ref{fig:Pu239_1}A shows that when increasing the number of ordinate $L$, the approximated eigenvalues from all three approaches converge to the ground truth value $k_{\text{eff}}=1$, as expected. Figure \ref{fig:Pu239_1}B shows the difference in approximated eigenvalues between the approaches. The full matrices iterative solver can reach to the accuracy around $10^{-8}$ compared to the one from generalized eigenvalue solver. The difference in eigenvalue between the full matrices and the TT approaches is around $10^{-7}$, which is below the TT truncated tolerance. Figure \ref{fig:Pu239_1}C shows the difference in the approximated eigenvectors. The full matrices and TT iterative approaches can reach to the accuracy around $10^{-4}$ compared to the eigenvector from generalized eigenvalue solver. The difference in eigenvector between the full grid and the TT approaches is around $10^{-8}$, which is again below the TT truncated tolerance.
Overall, this result shows that the TT approach produce very good approximations for both eigenvalues and eigenvectors.

\begin{figure}
    \centering
    \includegraphics[width=0.6\textwidth]{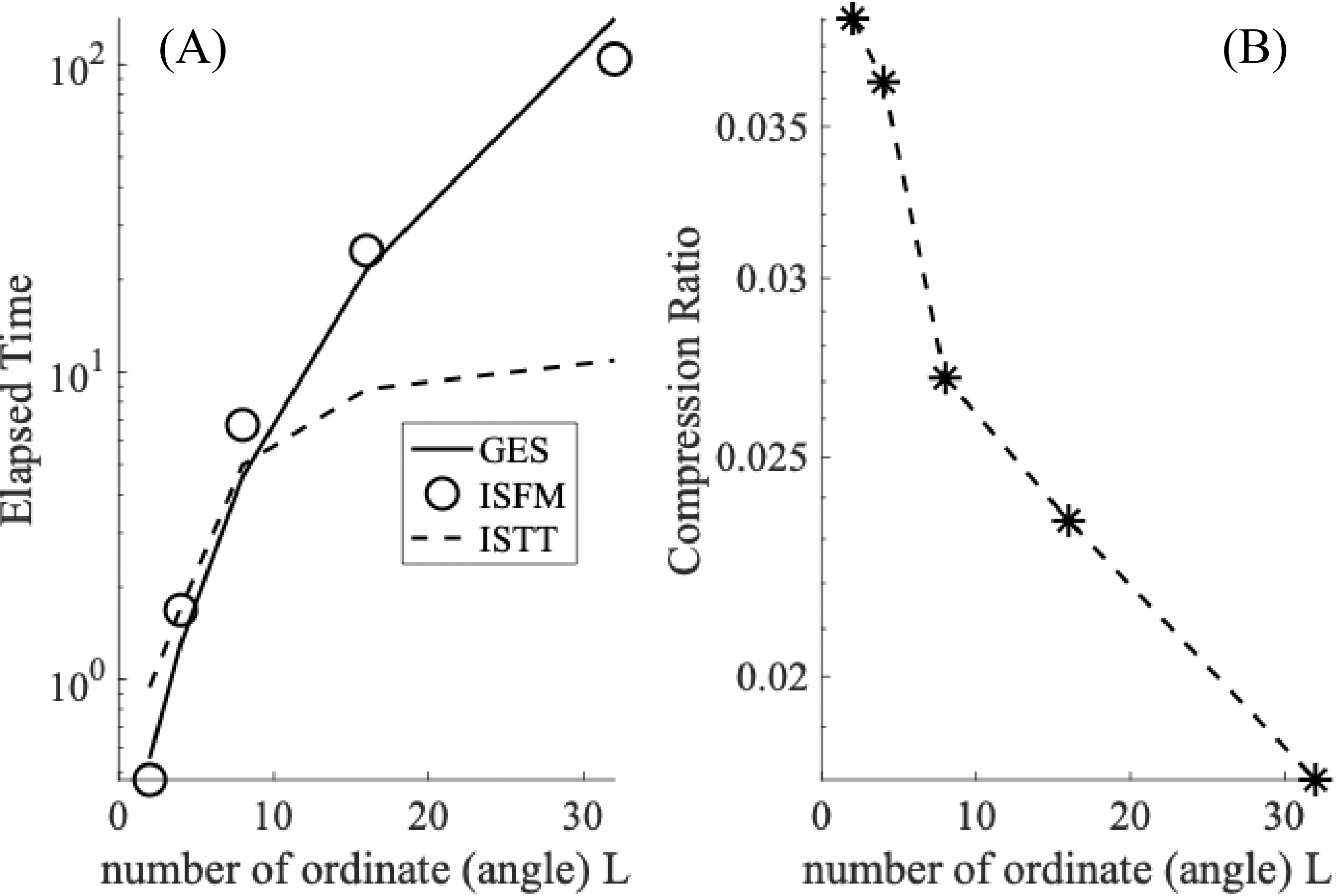}
    \caption{(\textit{Panel A}) - Comparison of the elapsed time between three approaches. The TT approach takes significantly less time compared to the others. At $L=32$, the TT approach is about 10 times faster (\textit{Panel B}) - Compression ratio of approximated eigenvectors in TT format shows that the TT approach already gains about 2 orders of magnitude in the storage cost.}
    \label{fig:Pu239_2}
\end{figure}

We also measure the compression of the TT format using the compression ratio, which is defined as follows:
\[\text{Compression ratio of } \Psi^{TT} = \dfrac{\# \text{ elements in TT format}}{\# \text{ elements in full grid tensor}}\]

Figure \ref{fig:Pu239_2} shows the elapsed time of each solver, and the compression ratio of the approximated eigenvector in the TT format.Figure \ref{fig:Pu239_2}A shows that the TT approach takes significantly less time compared to the other solvers. For example, for the problem size with $L=32$, the TT approach is about 10 times faster than the others. Figure \ref{fig:Pu239_2}B shows the compression ratios around $10^{-2}$, indicating that for these small problems, the TT format already gains about 2 orders of magnitude in term of the storage cost.

\subsubsection{Benchmark problems}
A variant of the critical assembly Jezebel was used to benchmark the tensor train approach for neutron transport. 
Jezebel was a tiny, nearly-spherical, nearly-bare (unreflected) experiment used from 1954-1955 to determine the critical mass of a homogeneous plutonium alloy \cite{favorite2021}. For this paper, a cube variant of the assembly was considered. The cube has dimensions \(x=(0,10)\), \(y=(0,10)\), \(z=(0,10)\). The material is a mixture of plutonium-239, plutonium-240, plutonium-241, gallium-69, and gallium-71. The precise composition of the material is listed in Table \ref{tab:benchmark}. The material composition is given in atomic density, defined as the number of atoms per cubic centimeter of the material.

\begin{table}[httb]
    \centering

\resizebox{\textwidth}{!}{
\begin{tabular}{|c|c|c|c|c|c|}
\hline
 &  \thead{Plutonium-239} & \thead{Plutonium-240} &  \thead{Plutonium-241} & \thead{Gallium-69} &  \thead{Gallium-71} \\
\hline
Atomic Density (\#/$cm^{3}$)& 3.7047e22 & 1.7512e21 & 1.1674e20 & 8.3603e20 & 5.3917e20 \\
\hline
\end{tabular}
}


    \caption{Critical Benchmark Material Composition}
    \label{tab:benchmark}
\end{table}

The cube dimensions were chosen such that the problem was slightly supercritical. Discretization of the problem used three different spatial grid sizes with $N_x = N_y = N_z = N \in \{128, 256, 512\}$ grid points per spatial dimension. The 128 quadrature points on the unit sphere are generated when using an $S_{8}$ quadrature \cite{walters1987}. Lastly, 256 energy groups are used. For deterministic neutron transport calculations at Los Alamos National Laboratory, 30 energy groups are traditionally used since the cost of the computation scales linearly with the number of energy groups. Using this discretization, the number of elements in the angular flux eigenvector $\Psi$ is approximately $6.8\times10^{10}$ elements (about 0.5 Terabytes (TB)) to $4.4\times10^{12}$ elements (about 4.4 TB). Both the k-effective and alpha-eigenvalue problems were considered.

The TT/QTT approach applied to these problems is implemented using functions from the  Oseledets's TT-Toolbox \cite{tt-toolbox} and is benchmarked on a Linux system with a 10-core i9 processor and 32GB RAM.

\subsubsection{Generation of Reference Values and PARTISN}

Reference values for the k-effective and alpha-eigenvalues were generated using PARTISN (PARallel TIme Dependent SN), the Los Alamos National Laboratory parallel time-dependent neutral particle $S{N}$ transport code package \cite{partisn2022}. The tensor train approach was compared to PARTISN. PARTISN numerically solves the neutron transport eigenvalue and source-driven problems in various geometries using a matrix-free solution method. The algorithm can briefly be described as follows: starting with an initial angular flux and eigenvalue, the scattering and fission sources are calculated for all cells, angles, and energy groups. Then, for all octants and energy groups, a known edge angular flux is used to "sweep" the problem by using known flux values to determine the neighboring cells' flux values. As mentioned earlier, unknown flux values are determined using the {\it diamond differencing} relationship. This solution of the angular flux iterate for all octants and energy groups is known as the "inner iteration." After obtaining the inner iteration fluxes, the scattering and fission sources are updated using the new angular flux, and the eigenvalue is updated. This is known as the "outer iteration." The process is repeated until the angular flux and eigenvalue converge to some tolerance. We note that the sweep of each energy group and octant is mathematically described as the solution of a lower triangular linear system using forward substitution. Throughout the algorithm, only the matrix-vector product of various quantities is required, allowing large problems to be solved without constructing the global linear system. We also note that the full neutron angular flux is not stored. Instead, the moments of the neutron angular flux are stored.
In most cases, only the first moment of the neutron angular flux, the neutron scalar flux, is retained for eigenvalue problems. This further reduces the amount of memory required. Parallelization is achieved by either a spatial or energy decomposition.

\subsubsection{The k-effective Eigenvalue Problem}
For these benchmark problems, the problem was decomposed over space in PARTISN. The performance of PARTISN is shown in Table \ref{tab:kPARTISN}. MPI ranks is the number of core that PARTISN used for the calculation on Snow, a LANL high-performance computing cluster.

\begin{table}[httb]
    \centering

\resizebox{\textwidth}{!}{
\begin{tabular}{|c|c|c|c|c|c|c|}
\hline
Grid Size &  \thead{Number of\\ Iterations} & \thead{Elapsed Time\\ (seconds)} &  \thead{Time per Iter \\ (seconds)} & \thead{MPI Ranks} &  \thead{Scalar Flux\\ Memory Storage (GB)} &  \thead{Eigenvalue} \\
\hline
128 &                    15 &        1131.323 &          75.42&                       72 &   4.29 &      1.0275823 \\
256 &                    15 &        5344.328 &         356.29&                      144 &   34.36 &      1.0276025 \\
512 &                    15 &        6061.987 &         404.13&                    1152 &   274.88 &      1.0276074 \\
\hline
\end{tabular}
}


    \caption{Performance of PARTISN on 3D benchmark k-effective eigenvalue problem at different grid sizes. MPI ranks are the numbers of cores that PARTISN used for the calculation.}
    \label{tab:kPARTISN}
\end{table}

Next, we applied the TT/QTT approach to the benchmark problems. To approximate how much speed-up the tensor train approach has gained, we calculate the speed-up factor defined as follows:
\begin{equation}
    \text{speed-up factor} = \dfrac{\text{PARTISN elapsed time}\times\text{Number of cores used by PARTISN}}{\text{TT/QTT elapsed time}\times\text{Number of cores used by TT/QTT}}
\end{equation}
\begin{table}[ht]
    \centering

\resizebox{\textwidth}{!}{
\begin{tabular}{|c|c|c|c|c|c|c|c|}
\hline
Grid Size &  \thead{Number of\\ Iterations} & \thead{Elapsed Time\\ (seconds)} & \thead{Speed-up\\ factor} & \thead{fullsize of $\ten{H}$ \\ (Zetabyte)} &  \thead{Compression Ratio \\ of $\ten{H}$ }  &  \thead{$\Psi^{TT}$ Memory \\ Storage (MB)} &  \thead{Eigenvalue\\ Error} \\
\hline
128 &                    16 &        37.58 &           216.75&   32 &   8.95e-18 &      1.03& 6.53e-5 \\
256 &                    22 &        57.07 &          1348.5& 2048 &   1.60e-19 &      0.98& 3.26e-5 \\
512 &                    18 &        92.32 &          7564.4& 131072 &   2.81e-21 &      1.15& 3.44e-5 \\
\hline
\end{tabular}
}


    \caption{Performance of Tensor Train approach on 3D benchmark k-effective eigenvalue problem at different grid sizes}
    \label{tab:kTT}
\end{table}

Table \ref{tab:kTT} shows that TT elapsed time is scaled at a slower rate when the grid size increases, increasing speed-up factors. For the largest problem at $512$, the TT approach is about 7500 times faster than PARTISN. Next, the compression rate of $\ten{H}$ shows the compression from yottabyte (YB) to megabyte (MB). The storage cost for $\Psi$ in QTT format is only around 1MB, compared to the full tensor storage cost of about 0.5 TB. Finally, the TT approach has the accuracy of $10^{-5}$ in eigenvalue compared to the PARTISN reference solution.

We also want to point out that in PARTISN, only the scalar flux is possibly retained with the storage cost of hundreds of GBs, and it is impossible to store the operators whose sizes are YB. While in the TT/QTT format, both the operators and the eigenvectors can be compressed and stored with the cost of MB. This will enable possibilities to use other methods to solve for multiple eigen-pairs at once.

\subsubsection{The Alpha-Eigenvalue Problem}
Next, PARTISN was used to generate the reference alpha-eigenvalue for the comparison. The performance of PARTISN is shown in Table \ref{tab:alphaPARTISN}.
\begin{table}[httb]
    \centering

\resizebox{\textwidth}{!}{
\begin{tabular}{|c|c|c|c|c|c|c|}
\hline
Grid Size &  \thead{Number of\\ Iterations} & \thead{Elapsed Time\\ (seconds)} &  \thead{Time per Iter \\ (seconds)} & \thead{MPI Ranks} &  \thead{Scalar Flux\\ Memory Storage (GB)} &  \thead{Eigenvalue} \\
\hline
128 &                    30 &        2475.821 &          82.53&                       72 &   4.29 &      9.2410746e-02 \\
256 &                    30 &        11637.28 &          378.91&                      144 &   34.36 &      9.2477425e-02 \\
512 &                    30 &        12783.57 &          426.12&                    1152 &   274.88 &      9.2493942e-02 \\
\hline
\end{tabular}
}


    \caption{Performance of PARTISN on 3D benchmark alpha-eigenvalue problems at different grid sizes}
    \label{tab:alphaPARTISN}
\end{table}

\begin{table}[httb]
    \centering


\resizebox{\textwidth}{!}{
\begin{tabular}{|c|c|c|c|c|c|c|c|}
\hline
Grid Size &  \thead{Number of\\ Iterations} & \thead{Elapsed Time\\ (seconds)} & \thead{Speed-up\\ factor} & \thead{fullsize of H \\ (Zetabyte)} &  \thead{Compression Ratio \\ of H}  &  \thead{$\Psi^{TT}$ Memory \\ Storage (MB)} &  \thead{Eigenvalue\\ Error} \\
\hline
128 &                    3 &        64.06 &           278.22&   32 &   8.95e-18 &      1.07& 1.13e-4 \\
256 &                    8 &        131.69 &          1272.6& 2048 &   1.60e-19 &      1.02& 2.89e-5 \\
512 &                    7 &        336.56 &          4375.6& 131072 &   2.81e-21 &      1.04& 4.44e-5 \\
\hline
\end{tabular}
}


    \caption{Performance of Tensor Train on 3D benchmark alpha-eigenvalue problem at different grid sizes}
    \label{tab:alphaTT}
\end{table}

The performance of the TT approach on the benchmark alpha-eigenvalue problems is shown in Table \ref{tab:alphaTT}. The result shows that the TT/QTT approach is significantly more efficient (about 4300 times faster for the largest problem). The increasing speed-up factor shows that the complexity of the TT approach scales at a lower rate compared to PARTISN. Moreover, the TT approach  
achieves the accuracy around $10^{-5}$ in approximating the alpha eigenvalues. The cost of storage is similar to the one from the k-effective eigenvalue problems, meaning while PARTISN needs about hundreds of GBs to store the scalar flux, the TT approach only uses about 1MB to store a whole eigenvector in TT format.


\section{Conclusions}
\label{sec:conclusions}

In this study, we utilize a low-rank tensor network method for
solution of ultra large time-independent integro-differential
Boltzmann Neutron Transport Equation.
We introduce a mixed Tensor Train (TT)/Quantized Tensor Train (QTT)
approach for numerical solution of 3D NTEs in Cartesian geometry.
We discretize the NTE based on $(a)$ diamond differencing; $(b)$
multigroup-in-energy; $(c)$ discrete ordinate collocation, which in
realistic cases leads to extremely large generalized eigenvalue
problems that requires matrix-free methods and large computer
clusters.
Further, we utilize this discretization and construct TT format of NTE
followed by QTT formating of the large TT-cores, which enable a
low-rank representation.
In the final QTT format of NTE we solve the tensorized generalized
eigenvalue problems using a fixed-point scheme and the TT-solver AMEn.
Comparing the full grid solutions, calculated by PARTISN, with the
solution of our TT/QTT method, we observed that the latter is
exceptionally efficient in terms of computational time and memory
usage, and significantly outperforms the full-grid version in
computational efficiency and storage requirements, achieving an
accuracy of $10^{-5}$.


\section*{Acknowledgements}
This work was supported by the Laboratory Directed Research and
Development (LDRD) program of Los Alamos National Laboratory under the
Grant 20230067DR, and in part by LANL Institutional Computing
Program. Los Alamos National Laboratory is operated by Triad National
Security, LLC, for the National Nuclear Security Administration of
U.S. Department of Energy (Contract No. 89233218CNA000001).

\bibliographystyle{plain}
\bibliography{main.bib}

\appendix
\section{}
\label{eqn:Appendix}
\subsection{The \textit{Multigroup-in-Energy} Approximation}
\label{subsubsec:multigroup-energy-approximation}

The discretization in the energy variable must be detailed enough to
capture the richness of neutron physics over a large energy domain.
For example, neutron resonances of various isotopes are regions where
the magnitude of a cross-section can vary rapidly over a small energy
range.
These nuclear physics effects must be carefully averaged in the energy
discretization process to ensure reaction probabilities remain
unchanged.
To discretize the energy variable E, we use the \textit{multigroup}
approximation; see, e.g., ~\cite{lewis_computational_1984}.
We restrict the energy $E$ to the finite interval
$\big[E_{\text{min}},E_{\text{max}}\big]$ that we partition into $G$
groups:
\begin{align*}
  E_{\text{max}} = E_{0} > E_{1} > \dots > E_{G} = E_{\text{min}}.
\end{align*}
Then, we average the eigenvalue equations
\eqref{eq:keff3DNeutronTransport} and
\eqref{eq:alpha3DNeutronTransport}, and we approximate the cross
sections by a flux-weighted average over each energy group $E_{g} < E<
E_{g-1}$.
We denote the neutron angular flux at the discretized energy group $g$
by $\psi_{g}$.

\subsection{The \textit{Discrete-Ordinates} Approximation}
\label{subsubsec:discrete-ordinates-approximation}

In nuclear reactors, the scattering of neutrons off materials like
water or graphite can be highly anisotropic, requiring many discrete
directions to resolve this behavior.
Therefore, we must evaluates the transport equation along a set of
discrete angular directions on the unit sphere that must be able to
capture the possible anisotropy of the neutron angular flux.
We discretize the angular variables using the \emph{discrete ordinates
method}~\cite{balshara2001OrdinateMethod}
To discretize the angular variable, we follow the process from
\cite{carlson1965transport} and consider a quadrature rule for
approximating integrals on the unit sphere $\mathcal{S}^{2}$:
\begin{align}
  \int_{\mathcal{S}^{2}}\diff\hat{\Omega}\,\psi(\hat{\Omega})
  \approx
  \sum_{\ell=1}^{L} w_{\ell} \psi(\hat{\Omega}_{\ell}) = \sum_{\ell=1}^{L} w_{\ell} \psi_{\ell},
  \label{eq:quadrature points}
\end{align}
where $\psi_{\ell}=\psi(\hat{\Omega}_{\ell})$ is the neutron angular
flux at direction $\ell$ and
$\hat{\Omega}_{\ell}\equiv(\mu_{\ell},\eta_{\ell},\xi_{\ell})$, for
$\ell$ ranging through $1$ and $L=2N^{2}$, $N$ being the number of
direction cosines.
We assume that the quadrature weights $w_{\ell}$ are normalized so that
$\sum_{\ell=1}^{L} w_{\ell} = 1$.
In this work, We assume a square quadrature: a quadrature set where each axis has the same quadrature points forming a square grid of unique ordinates. Other types of quadrature can be considered if desired.
\begin{align*}
  \begin{split}
    -1 < \mu_{1}  < \dots < \mu_{N/2}  < 0 < \mu_{N/2+1}  < \dots < \mu_{N}  < 1, \mu_{N+1-n}  = - \mu_{n}, \\[0.5em]
    -1 < \eta_{1} < \dots < \eta_{N/2} < 0 < \eta_{N/2+1} < \dots < \eta_{N} < 1, \eta_{N+1-n} = - \eta_{n}.
  \end{split}
\end{align*}
We note that since $\hat{\Omega}_{\ell}\in\mathcal{S}^{2}$ for all
$\ell$, $\xi_{l}$
and $\xi_{\ell}=\sqrt{1-\mu_{\ell}^{2}-\eta_{\ell}^{2}}$ 
, is given by
\begin{align*}	
  \xi_{\ell} = \sqrt{ 1 - \mu_{\ell}^{2} - \eta_{\ell}^{2} }
\end{align*}
We collect the ordinates $\mu_{\ell}$ in the arrays
$\hat{\mu}_{-}=\big[\mu_{1},\ldots,\mu_{L/2}\big]^{T}$ and
$\hat{\mu}_{+}=\big[\mu_{L/2+1},\ldots,\mu_{L}\big]^{T}$, and we
similarly define the arrays $\hat{\eta}_{-}$, $\hat{\eta}_{+}$,
$\hat{\xi}_{-}$, and $\hat{\xi}_{+}$.
The ordering of the ordinates is arbitrary. The various combinations of ordinates specify faces of a three-dimensional surface and the ordering is chosen such that boundary conditions are simple to specify. For example, in reactor calculations, reactor fuel elements are symmetric and these lines of symmetry can be be used to reduce the size of the problem. In this case, the ordering of the ordinates prioritizes the vacuum boundary condition faces first and then the reflective faces. In this work we choose the ordinates ordering such that
\begin{equation}
  \left\{\,
  \begin{array}{l}
    \big[ \hat{\mu}_{-}, \hat{\eta}_{-}, \hat{\xi}_{-} \big] \\[0.5em]
    \big[ \hat{\mu}_{+}, \hat{\eta}_{-}, \hat{\xi}_{-} \big] \\[0.5em]
    \big[ \hat{\mu}_{-}, \hat{\eta}_{+}, \hat{\xi}_{-} \big] \\[0.5em]
    \big[ \hat{\mu}_{+}, \hat{\eta}_{+}, \hat{\xi}_{-} \big] \\[0.5em]
    \big[ \hat{\mu}_{-}, \hat{\eta}_{-}, \hat{\xi}_{+} \big] \\[0.5em]
    \big[ \hat{\mu}_{+}, \hat{\eta}_{-}, \hat{\xi}_{+} \big] \\[0.5em]
    \big[ \hat{\mu}_{-}, \hat{\eta}_{+}, \hat{\xi}_{+} \big] \\[0.5em]
    \big[ \hat{\mu}_{+}, \hat{\eta}_{+}, \hat{\xi}_{+} \big]
  \end{array}
  \,\right\}
  \in \mathbb{R}^{L \times 3},
  \label{eqn:hat-mu}
\end{equation}
since we are dealing exclusively with vacuum boundary conditions.
\subsubsection{The \textit{Diamond-Difference} Approximation}
\label{subsubsec:diamond-difference-approximation}

The space discretization must be of the order of the mean free path of
a neutron before interacting, i.e., millimeters in a system like a
nuclear reactor which measures in meters.
We perform the discretization on the space independente variables $x$,
$y$, and $z$ using the \emph{diamond differencing method}.

In the spatial dimension, we introduce a univariate grid partition
over each problem dimension.
We consider the grid stepsizes $\Delta x_i$, $\Delta y_j$, and $\Delta z_k$
to partition the real, bounded intervals $\big[a_x,b_x\big]$,
$\big[a_y,b_y\big]$, and $\big[a_z,b_z\big]$ into $M$, $J$, and $K$
cells, respectively, so that
\begin{align*}
  \begin{array}{rll}
    a_{x}&\!\!\!\!\equiv x_{0} < \ldots < x_{i-1} < x_{i} < \ldots < x_{M} &\!\!\!\!\!\equiv b_{x}, \\[0.25em]
    a_{y}&\!\!\!\!\equiv y_{0} < \ldots < y_{j-1} < y_{j} < \ldots < y_{J} &\!\!\!\!\!\equiv b_{y}, \\[0.25em]
    a_{z}&\!\!\!\!\equiv z_{0} < \ldots < z_{k-1} < z_{k} < \ldots < z_{K} &\!\!\!\!\!\equiv b_{z},
  \end{array}
\end{align*}
and $\Delta x_{i} = x_{i} - x_{i-1}$, $\Delta y_{j} = y_{j} - y_{j-1}$, and 
$\Delta z_{k} = z_{k} - z_{k-1}$.
We refer to the grid nodes $x_{i}$, $y_{j}$, $z_{k}$ as
``\emph{edges}'' and we call ``\emph{edge values}'' the corresponding
function values.

\medskip
Since cross-sections are constant in a cell, the cell-centered angular
flux must be expressed in terms of edge angular flux values.
The cell-centered angular flux is located at $i-1/2$ and the {\it
  diamond difference} approximation~\cite{lewis_computational_1984}
gives the cell-centered angular fluxes at the edges labeled by the
index sets $(k,j,i-1/2)$, $(k,j-1/2,i)$, and $(k-1/2,j,i)$ as
\begin{align*}
  \psi_{g,\ell,k,j,i-1/2} = \frac{1}{2} \big ( \psi_{g,\ell,k,j,i} + \psi_{g,\ell,k,j,i-1} \big ), \\
  \psi_{g,\ell,k,j-1/2,i} = \frac{1}{2} \big ( \psi_{g,\ell,k,j,i} + \psi_{g,\ell,k,j-1,i} \big ), \\
  \psi_{g,\ell,k-1/2,j,i} = \frac{1}{2} \big ( \psi_{g,\ell,k,j,i} + \psi_{g,\ell,k-1,j,i} \big ).
\end{align*}
\subsubsection{Discretization of the boundary conditions}
Boundary conditions are defined for all faces of the three-dimensional problem. An incoming neutron angular flux can be imposed on a face of the problem. For example, an boundary flux on the top face of the three-dimensional cube would be expressed as
\begin{equation}
 \psi_{g,\ell,K,J,M} = \Psi_{\text{top}} \quad\text{ for } \mu_{\ell} < 0, \eta_{\ell} < 0, \xi_{\ell} < 0,
\end{equation}
since we have defined $(\mu_{\ell} < 0, \eta_{\ell} < 0, \xi_{\ell} < 0)$ to be the inward direction, and $K,J,M$ specifies that this incoming angular flux is for all cells in the top face of the problem. For eigenvalue problems, the boundary conditions are taken to be vacuum boundary conditions (no incoming angular flux) and the discretized boundary conditions for each face are given by

\begin{equation}
  \begin{split}
    \psi_{g,\ell,K,J,M} &= 0 \quad\text{ for } \mu_{\ell} < 0, \eta_{\ell} < 0, \xi_{\ell} < 0 \\
    \psi_{g,\ell,K,J,0} &= 0 \quad\text{ for } \mu_{\ell} > 0, \eta_{\ell} < 0, \xi_{\ell} < 0 \\[0.25em]
    \psi_{g,\ell,K,0,M} &= 0 \quad\text{ for } \mu_{\ell} < 0, \eta_{\ell} > 0, \xi_{\ell} < 0 \\[0.25em]
    \psi_{g,\ell,K,0,0} &= 0 \quad\text{ for } \mu_{\ell} > 0, \eta_{\ell} > 0, \xi_{\ell} < 0 \\[0.25em]
    \psi_{g,\ell,0,J,M} &= 0 \quad\text{ for } \mu_{\ell} < 0, \eta_{\ell} < 0, \xi_{\ell} > 0 \\[0.25em]
    \psi_{g,\ell,0,J,0} &= 0 \quad\text{ for } \mu_{\ell} > 0, \eta_{\ell} < 0, \xi_{\ell} > 0 \\[0.25em]
    \psi_{g,\ell,0,0,M} &= 0 \quad\text{ for } \mu_{\ell} < 0, \eta_{\ell} > 0, \xi_{\ell} > 0 \\[0.25em]
    \psi_{g,\ell,0,0,0} &= 0 \quad\text{ for } \mu_{\ell} > 0, \eta_{\ell} > 0, \xi_{\ell} > 0.
  \end{split}
  \label{eq:discreteBCs}
\end{equation}
The set of equations defined by Equations \ref{eq:Discretized_kNTE} and \ref{eq:Discretized_alphaNTE} have $GLKJM$ unknowns and $GLMK + GLJM + GLJK + GLM + GLJ + GLK + GL$ boundary equations.:




\subsection{Matricization of the Discrete Neutron Transport Eigenvalue Equations}
To write Equations \ref{eq:Discretized_kNTE} and \ref{eq:Discretized_alphaNTE} in matrix form, we define the angular flux vector for a single energy group $g$ and direction $\ell$ as
\begin{multline}
\Psi_{g, \ell} \equiv 
\begin{pmatrix}
\psi_{g,\ell,0} \\
\vdots \\
\psi_{g,\ell,K}
\end{pmatrix} \in \mathbb{R}^{(K+1)(J+1)(M+1)}, \\
\Psi_{g, \ell, k} \equiv 
\begin{pmatrix}
\psi_{g,\ell,k,0} \\
\vdots \\
\psi_{g,\ell,k,J}
\end{pmatrix} \in \mathbb{R}^{(J+1)(M+1)}, \\
\Psi_{g,\ell,k,j} \equiv 
\begin{pmatrix}
\psi_{g,\ell,k,j,0} \\
\vdots \\
\psi_{g,\ell,k,j,M}
\end{pmatrix} \in \mathbb{R}^{(M+1)}.
\label{eq:aflux}
\end{multline}
The angular flux vector for an energy group $g$ is defined as
\begin{equation}
\Psi_{g} \equiv 
\begin{pmatrix}
\Psi_{g,1} \\
\vdots \\
\Psi_{g,L}
\end{pmatrix} \in \mathbb{R}^{L(K+1)(J+1)(M+1)}.
\end{equation}
The full angular flux vector is then defined as
\begin{equation}
\Psi \equiv
   \begin{pmatrix}
     \Psi_{1} \\
     \Psi_{2} \\
     \vdots \\
     \Psi_{G}
   \end{pmatrix} \in \mathbb{R}^{GL(K+1)(J+1)(M+1)}.
\end{equation}

\subsubsection{Differencing and Interpolation Matrices}
To write the matrix form of the diamond difference discretized derivative operators $\mu \partial/\partial x + \eta \partial/\partial y + \xi \partial/\partial z$, we define the matrices
\begin{equation}
\begin{split}
\Delta x \equiv &\text{diag}(\Delta x_{1}, \dots, \Delta x_{M}) \in \mathbb{R}^{M \times M}, \\
\Delta y \equiv &\text{diag}(\Delta y_{1}, \dots, \Delta y_{J}) \in \mathbb{R}^{J \times J}, \\
\Delta z \equiv &\text{diag}(\Delta z_{1}, \dots, \Delta z_{K}) \in \mathbb{R}^{K \times K}.
\end{split}
\end{equation}
The differencing matrix is defined as
\begin{equation}
D_{x} \equiv
\setlength\arraycolsep{2pt}
\begin{pmatrix}
-1 & 1 & & \\
& \ddots & \ddots & \\
& & -1 & 1
\end{pmatrix} \in \mathbb{R}^{M \times (M+1)},
\end{equation}
where the matrices $D_{y}$ and $D_{z}$ are defined similarly. Since each derivative term requires an average of cell-centered angular fluxes, we define the matrix $S_{x}$ as
\begin{equation}
S_{x} = \frac{1}{2}
\setlength\arraycolsep{2pt}
\begin{pmatrix}
1 & 1 & & \\
& \ddots & \ddots & \\
& & 1 & 1
\end{pmatrix} \in \mathbb{R}^{M \times (M+1)},
\end{equation}
where $S_{y}$ and $S_{z}$ are defined similarly. The matrices $S_{x}$, $S_{y}$, and $S_{z}$ interpolate cell-centered vectors into zone-centered vectors by averaging the edge angular flux values.

\subsubsection{Cell/Edge Transformation Matrices}
We now define the matrices $Z$ and $Z_{b}$ as
\begin{equation}
	Z \equiv \begin{pmatrix}
			I_{MJK} \\
			0
		     \end{pmatrix} \in \mathbb{R}^{(K+1)(J+1)(M+1) \times KJM},
\end{equation}
\begin{equation} 
	Z_{b} \equiv \begin{pmatrix}
			0 \\
			I_{(K+1)(J+1)(M+1)-KJM)}
		      \end{pmatrix} \in \mathbb{R}^{(K+1)(J+1)(M+1) \times (K+1)(J+1)(M+1)-KJM}.
\end{equation}
The matrices $Z$ and $Z_{b}$ transform cell-centered vectors to edge vectors and vice versa.

\subsubsection{The Discrete Ordinates Matrices}
We define the angular quadrature point matrices as
\begin{equation}
\begin{split}
\hat{\mu} &= \text{diag}( \hat{\mu}_{-}, \hat{\mu}_{+}, \hat{\mu}_{-}, \hat{\mu}_{+} ) \in \mathbb{R}^{L \times L}, \\
\hat{\eta} &= \text{diag}( \hat{\eta}_{-}, \hat{\eta}_{-}, \hat{\eta}_{+}, \hat{\eta}_{+} ) \in \mathbb{R}^{L \times L}, \\
\hat{\xi} &= \text{diag}( \hat{\xi}_{-}, \hat{\xi}_{-}, \hat{\xi}_{+}, \hat{\xi}_{+} ) \in \mathbb{R}^{L \times L}. \\
\end{split}
\label{eqn:hatmu}
\end{equation}

Discretized representations of the angular flux moment operators must be defined. These operators operate on zone-centered vectors and are easily seen to be given by $KJM \times LKJM$ size matrices
\begin{equation}
	L_{n,m} \equiv (l_{n,m}W) \otimes I_{MJK},
\end{equation}
where
\begin{equation}
	l_{n,m} \equiv \bigg (Y_{n}^{m}(\hat{\Omega}_{1}), Y_{n}^{m}(\hat{\Omega}_{2}), \dots, Y_{n}^{m}(\hat{\Omega}_{L}) \bigg ),
\end{equation}
and
\begin{equation}
	W \equiv \text{diag}(w_{1}, w_{2}, \dots, w_{L}).
\end{equation}
If the vector $\Psi$ approximates $\psi(\vec{r}, \hat{\Omega})$, then $L_{n,m} \Psi$ approximates the (n,m)$^{\text{th}}$ moment of $\psi(\vec{r},\hat{\Omega})$, $\phi_{n,m}(\vec{r})$. Similarly, we define $LMJK \times MJK$ size matrices
\begin{equation}
	L^{+}_{n,m} \equiv l^{T}_{n,m} \otimes I_{MJK}.
\end{equation} 
If a vector $\Phi$ approximates $\phi(\vec{r})$, then $L^{+}_{n,m} \Phi$ approximates $Y_{n}^{m}(\hat{\Omega}) \phi(\vec{r})$. We define the grouped matrices $L_{n}$ and $L^{+}_{n}$, where
\begin{equation}
	L_{n} = \begin{pmatrix}
			L_{n,-n} \\
			\vdots \\
			L_{n,n}
		    \end{pmatrix} \text{ and }
	L^{+}_{n} = \begin{pmatrix}
		    	L^{+}_{n,-n}, \dots, L^{+}_{n,n},
		    \end{pmatrix}
\end{equation}
and the further grouped block matrices
\begin{equation}
	L^{N} = \begin{pmatrix}
			L_{0} \\
			\vdots \\
			L_{N}
		    \end{pmatrix} \text{ and }
	L^{N,+} = \begin{pmatrix}
				L^{+}_{0}, \dots, L^{+}_{N}
			\end{pmatrix},
\end{equation}
where $N = N_{s}$, the number of terms in the scattering kernel. The scattering cross section in Eq.~\ref{eq:keff3DNeutronTransport} is usually expanded in spherical harmonics up to some order $N_{s}$\cite{bell_nuclear_1970}. For this reason, it is assumed that the symmetric quadrature rule is such that the spherical harmonics of order $N_{s}$ and less satisfy \cite{brown_linear_1995}
\begin{equation}
	\sum_{\ell = 1}^{L} Y_{n}^{m}(\hat{\Omega}_{\ell}) Y_{n'}^{m'}(\hat{\Omega}_{\ell}) = \delta_{n,n'} \delta_{m,m'}, \text{ for all } 0 \leq n, n' \leq N_{s}, \abs{m} \leq n, \abs{m'} \leq n'.
\end{equation}
In this work, we assume all scattering is isotropic ($N_{s}=0$) for simplicity.

\subsubsection{Matrix Representations of the Spatial Derivatives and Total Cross Section}
The matrix form representation of  $\mu \, \partial/\partial x$, $H_{\mu}$, can be written as
\begin{equation}
        H_{\mu} = I_{G} \otimes \hat{\mu} \otimes Z ( S_{z} \otimes S_{y} \otimes \Delta x^{-1}D_{x} ),
\end{equation}
where $I_{G}$ is the identity matrix sized to the number of energy groups. The matrix $H_{\mu}$ is a square matrix with size $GL(K+1)(J+1)(M+1)\times GL(K+1)(J+1)(M+1)$. The other two derivative matrix terms, $H_{\eta}$ and $H_{\xi}$ can be written as
\begin{equation}
        H_{\eta} = I_{G} \otimes \hat{\eta} \otimes Z ( S_{z} \otimes \Delta y^{-1}D_{y} \otimes S_{x} ),
\end{equation}
and
\begin{equation}
        H_{\xi} = I_{G} \otimes \hat{\xi} \otimes Z ( \Delta z^{-1}D_{z} \otimes S_{y} \otimes S_{x} ).
\end{equation}
We define the total cross section matrices for energy group $g$ over all cells as
\begin{equation}
	\Sigma_{g} \equiv \text{diag}(\sigma_{g,111},\dots,\sigma_{g,KJM}) \in \mathbb{R}^{KJM \times KJM}.
\end{equation}
The total cross section matrix for all energy groups is then
\begin{equation}
\Sigma = ( I_{G} \otimes I_{L} \otimes Z ) \bigg ( I_{L} \otimes \text{diag} ( \Sigma_{1}, \dots, \Sigma_{G} ) \bigg ) \big ( I_{G} \otimes I_{L} \otimes S \big )
\end{equation}
where $S = S_{z} \otimes S_{y} \otimes S_{x}$.

To apply the boundary conditions shown in Eq.~\ref{eq:discreteBCs}, we define the matrices $E_{kji}$ that pick out the correct elements of $\Psi$ for some $\hat{\Omega}_{\ell}$ as done in \cite{brown_linear_1995}. There are eight different $E_{kji}$ matrices in three-dimensions with $k=0$ or $K$, $j=0$ or $J$, and $i=0$ or $M$. For vacuum boundary conditions, we have
\begin{equation}
    E_{kji}\Psi = 0.
\end{equation}
For ordinate $\hat{\Omega}_{\ell} > 0$ $(\mu_{\ell}, \eta_{\ell}, \xi_{\ell} > 0)$, the boundary matrix $E_{000}$ is
\begin{equation}
E_{000} =
\begin{pmatrix}
e_{0K}^{T} \otimes I_{J+1} \otimes I_{M+1} \\
(0,I_{K}) \otimes e_{0J}^{T} \otimes I_{M+1} \\
(0,I_{K}) \otimes (0,I_{J}) \otimes e_{0M}^{T}
\end{pmatrix},
\end{equation}
where $e_{0M}$ is the basis vector sized $M+1$ with one as the first element and zero elsewhere. The other basic vectors are defined similarly. The boundary condition matrix $B$ is then defined as
\begin{equation}
    B = ( I_{G} \otimes I_{L} \otimes Z_{b} ) \bigg ( I_{G} \otimes \text{diag}( E_{KJM}, E_{KJM}, \dots, E_{000}, E_{000} ) \bigg ),
\end{equation}
where the ordering of the matrices $E_{kji}$ is determined by the signs of the quadrature points $(\mu_{\ell}, \eta_{\ell}, \xi_{\ell})$.

\subsubsection{Matrix Representation of the Scattering Cross Section}
The scattering cross section matrix is defined by letting
\begin{equation}
	\Sigma_{s,g,g'} =  \text{diag}(\sigma_{s,g,g',111}, \dots, \sigma_{s,g,g',KJM}).
\end{equation}
The scattering cross section for all energy groups is then
\begin{equation}
	\Sigma_{s} \equiv I_{L} \otimes \begin{pmatrix}
							\Sigma_{s,11} & \dots & \Sigma_{s,1G} \\
							\vdots & \ddots & \vdots \\
							\Sigma_{s,G1} & \dots & \Sigma_{s,GG}
						\end{pmatrix}.
\end{equation}
The full scattering operator $\mathbb{S}$ is then
\begin{equation}
    \mathbb{S} = \big ( I_{G} \otimes I_{L} \otimes Z \big ) \big ( I_{G} \otimes L^{0,+} \big ) \Sigma_{s} \big ( I_{G} \otimes L^{0} \big ) \big ( I_{G} \otimes I_{L} \otimes S \big ).
\end{equation}

\subsubsection{Matrix Representation of the Fission Cross Section}
The fission cross section matrix is defined similarly
\begin{equation}
	\Sigma_{f,g,g'} \equiv \text{diag}(\chi_{g'g} \nu \sigma_{f,g',111}, \dots, \chi_{g'g} \nu \sigma_{f,g',KJM}),
\end{equation}
where the fission cross section for all energy groups is then
\begin{equation}
	\Sigma_{f} \equiv I_{L} \otimes  \begin{pmatrix}
							\Sigma_{f,11} & \dots & \Sigma_{f,1G} \\
							\vdots & \ddots & \vdots \\
							\Sigma_{f,G1} & \dots & \Sigma_{f,GG}
						\end{pmatrix}.
\end{equation}

The fission operator $\mathbb{F}$ is then
\begin{equation}
    \mathbb{F} = \big ( I_{G} \otimes I_{L} \otimes Z \big ) \big ( I_{G} \otimes L^{0,+} \big ) \Sigma_{f} \big ( I_{G} \otimes L^{0} \big ) \big ( I_{G} \otimes I_{L} \otimes S \big ).
\end{equation}

\subsubsection{Matrix Representation of the Inverse Neutron Group Velocity}

The inverse velocity cross section matrices for energy group $g$ over all cells is given by
\begin{equation}
	V_{g}^{-1} \equiv \frac{1}{v_{g}} I_{KJM} \in \mathbb{R}^{KJM \times KJM},
\end{equation}
where $I$ is identity matrix with size $KJM$. The inverse velocity matrix for all energy groups is then
\begin{equation}
\mathbb{V}^{-1} = ( I_{G} \otimes I_{L} \otimes Z ) \bigg ( I_{L} \otimes \text{diag} ( V_{1}^{-1}, \dots, V_{G}^{-1} ) \bigg ) \big ( I_{G} \otimes I_{L} \otimes S \big ).
\end{equation}


\section{} 
\label{Appendix_B}
\subsection{Notation, basic definitions, and operations with tensors}
Let $d$ be a positive integer.
A $d$-dimensional tensor $\ten{A}\in\mathbb{R}^{n_1\times
  n_2\times\ldots\times n_d}$ is a multi-dimensional array with $d$
indices and $n_k$ elements in the $k$-th direction, $k=1,2,\ldots,d$
being the dimensional index.
We say that the number of dimensions $d$ is the \emph{order of the
tensor}.
As usual, we refer to one-dimensional tensors as \emph{vectors}, and
two-dimensional tensors as \emph{matrices}.
We denote the tensors using uppercase, calligraphic fonts, e.g.,
$\ten{A}$; the matrices with bold, uppercase fonts, e.g., $\mat{A}$;
the vectors with bold, lowercase fonts, e.g., $\mat{a}$.
To denote tensor's, matrix's, and vector's components, we use both the
subscripted notation, e.g., $\ten{A}=(\mathcal{A}_{ijk})$,
$\mat{a}=(a_i)$, $\mat{A}=(A_{ij})$, and the 
MATLAB$^{\copyright}$~\cite{MATLAB:2019} notation, e.g.,
\begin{equation}
  \ten{A}:= \big( \ten{A}(i_1,i_2,\dots,i_d) \Big),
  \quad i_k = 1,\dots,n_k,
  \quad k   = 1,\dots, d;
  \label{eq:tns:notation}
\end{equation}
$\mat{A}=\big(\mat{A}(i,j)\big)\in\mathbb{R}^{m\times n}$ and
$i=1,\ldots,m$, $j=1,\ldots,n$;
$\mat{a}=\big(\mat{a}(i)\big)\in\mathbb{R}^{n}$ and $i=1,\ldots,n$.

We form a tensor subarray by fixing one or more of its indices.
For example, the \emph{tensor fibers} (the higher-order analog of
matrix rows or columns) are defined by fixing all but one of the
tensor indices, while the \emph{tensor slices} are two-dimensional
sections, defined by fixing all but two of the tensor indices.
For example, still using a MATLAB-like notation, $\ten{A}(i_1,i_2,:)$ and $\ten{A}(i_1,:,:)$ respectively denote the fiber along the third direction and the slice for any fixed index value $i_1$.
 

\subsubsection{Kronecker product}
The Kronecker product $\bigotimes$ of matrix
$\mat{A}=(a_{ij})\in\mathbb{R}^{m_A \times n_A}$ and matrix
$\mat{B}=(b_{ij})\in\mathbb{R}^{m_B \times n_B}$ is the matrix
$\mat{A}\otimes \mat{B}$ of size
$N_{\mat{A}\otimes\mat{B}}=(m_Am_B)\times(n_An_B)$ defined as:
\begin{equation}
  \mat{A}\otimes\mat{B} = 
  \begin{bmatrix}
    a_{11}\mat{B} &a_{12}\mat{B} &\cdots & a_{1n_A}\mat{B}\\
    a_{21}\mat{B} &a_{22}\mat{B} &\cdots & a_{2n_A}\mat{B}\\
    \vdots & \vdots & \ddots & \vdots \\
    a_{m_A1}\mat{B} &a_{m_A2}\mat{B} &\cdots & a_{In_A}\mat{B}\\
  \end{bmatrix}.
  \label{def:kronecker}
\end{equation}
Equivalently, it holds that
$\big(\mat{A}\otimes\mat{B}\big)_{ij}=a_{i_Aj_A}b_{i_Bj_B}$, where
$i=i_B+(i_A-1)m_B$, $j=j_B+(j_A-1)m_B$, with $i_A=1,\ldots,m_A$,
$j_A=1,\ldots,n_A$, $i_B=1,\ldots,m_B$, and $j_B=1,\ldots,n_B$.

\subsubsection{Tensor product}
 There is a relation between Kronecker product and tensor product: Kronecker product is a \textit{particular} bilinear map on a pair of vector spaces consisting of matrices of a given dimensions (it requires a choice of basis), while the tensor product is a \textit{universal} bilinear map on a pair of vector spaces of any sort (i.e., it is more general). 

Here we define the tensor product of two vectors $\mat{a}=(a_i)\in\mathbb{R}^{n_A}$
and $\mat{b}=(b_i)\in\mathbb{R}^{n_B}$, which produces the matrix
$\big(\mat{a}\tenprod\mat{b}\big)$ of size
$N_{\mat{a}\tenprod\mat{b}}=n_A\times n_B$ defined as:
\begin{equation}
  \big(\mat{a}\tenprod\mat{b}\big)_{ij} = a_{i}b_{j}
  \qquad
  i=1,2,\ldots,n_A,\,
  j=1,2,\ldots,n_B.
  \label{def:outer:product:vectors}
\end{equation}
Note that
$\mat{a}\tenprod\mat{b}=\mat{a}\otimes\mat{b}^T$.
 Similarly, the tensor product of matrix
$\mat{A}=(a_{ij})\in\mathbb{R}^{m_A \times n_A}$ and matrix
$\mat{B}=(b_{kl})\in\mathbb{R}^{m_B \times n_B}$ produces the
four-dimensional tensor of size $N_{\mat{A}\tenprod\mat{B}}=m_A \times
n_A \times m_B \times n_B$, with elements:
\begin{equation}
  \big(\mat{A}\tenprod\mat{B}\big)_{ijkl} = a_{ij}b_{kl},
  \label{def:outer_product:matrices}
\end{equation}
for
$i=1,2,\ldots,m_A$,
$j=1,2,\ldots,n_A$,
$k=1,2,\ldots,m_B$,
$l=1,2,\ldots,n_B$.

\subsubsection{Contraction of a tensor with a vector}
\label{subsubsec:contraction:product}
Consider the tensor
$\ten{X}\in\mathbb{R}^{n_1\times n_2\times\dots\times n_d}$
and the vector $\mat{v}\in\mathbb{R}^{n_k}$ for some $1\leq k\leq d$.
The \emph{$k$-th tensor-vector contraction} of $\ten{X}$ with
$\mat{v}$ is the summation over the $k$-th index of the tensor
elements weighted by the vector components:
\begin{equation}
  \big(\ten{X}\Bar{\times}_k\mat{v}\big)(i_1,i_2,\ldots,i_{k-1},i_{k+1},\ldots i_d) =
  \sum_{i_k=1}^{n_k}
  \ten{X}(i_1,i_2,\ldots,i_{k-1},i_k,i_{k+1},\ldots i_d)
  \mat{v}(i_k).
\end{equation}
Tensor $\ten{X}\Bar{\times}_k\mat{v}$ is a $(d-1)$-dimensional array
of size
$N_{\ten{X}\Bar{\times}_k\mat{v}} = n_1\times n_2\times\ldots\times n_{k-1}\times n_{k+1}\times\ldots\times\
n_d$.

\subsubsection{The $n$-mode product}
\label{subsubsec:n-mode:product}
Consider the tensor
$\ten{X}\in\mathbb{R}^{n_1\times n_2\times\dots\times n_d}$
and the matrix $\mat{U}\in\mathbb{R}^{n_U\times n_n}$.
The \emph{$n$-mode product} between $\ten{X}$ and $\mat{U}$ is the
contraction along the $n$-th direction given by
\begin{multline*}
  \big(\ten{X}\times_n\mat{U}\big)(i_1,\ldots,i_{n-1},\ell,i_{n+1},\ldots,i_d) =
  \\[0.5em]
  = \sum^{n_n}_{j=1} \ten{X}(i_1,\ldots,i_{n-1},j,i_{n+1},\ldots,i_{d})\mat{U}(\ell,j)
  \quad
  \forall \ell=1,2,\ldots,n_U.
\end{multline*}
Tensor $\ten{X}\times_n\mat{U}$ has the same dimension $d$ of $\ten{X}$,
but size $N_{\ten{X}\times_n\mat{U}} = n_1\times n_2\times\ldots\times
n_{n-1}\times n_U\times n_{n+1}\times n_d$ instead of $N_{\ten{X}} =
n_1\times n_2\times\ldots n_d$, which is the size of $\ten{X}$.

\subsection{Differential Operators in \TT{format}.}
\label{APP:Diff_Operators}
Applying, for example, the first-order accurate, forward difference formula to
approximate the differentiation of $f(x_1,x_2,x_3,x_4)$ through the values of 
tensor $\ten{F}$ along the direction of the independent variable $x_2$ yields:
\begin{align*}
  \left(\dfrac{\partial f}{\partial x_2}\right)_{i_1,i_2,i_3,i_4}
  = \dfrac{\ten{F}(i_1,i_2+1,i_3,i_4)-\ten{F}(i_1,i_2,i_3,i_4)}{\Delta x_2}
  + \mathcal{O}\big(\Delta x_2\big),
\end{align*}
for $i_2=1,2,\ldots,(n_2-1)$, and where $\Delta x_2$ is the grid
step-size along the direction of the independent variable $x_2$.
When working with grid functions, all the discrete analogs of the
differential operators (i.e., gradient, curl, divergence, Laplacian,
etc.) must be expressed in \TT{format}.
Because of the discrete separation of variables provided by the
\TT{format}, the forward difference scheme approximating
$\partial{f}/\partial{x_2}$ is acting only on the index $i_2$ of the
second \TT{core} $\ten{G}_2(:,i_2,:)$ of $\ten{G}^{TT}$ and can be
computed directly in the \TT{format} as follows:
\begin{multline*}
  \left(\dfrac{\partial f}{\partial x_2}\right)^{TT}_{i_1,i_2,i_3,i_4}
  =
  \sum_{\alpha_{1},\alpha_{2},\alpha_{3}=1}^{r_1,r_2,r_3}
  \ten{G}_1(1,i_1,\alpha_1)
  \dfrac{\ten{G}_2(\alpha_{1},i_2+1,\alpha_2) - \ten{G}_2(\alpha_{1},i_2,\alpha_2)}{\Delta{x_2}}
  \ldots\\ \ldots
  \ten{G}_3(\alpha_{2},i_3,\alpha_3)\ten{G}_4(\alpha_{3},i_4,1)
  + \varepsilon,
\end{multline*}
where the ``dots'' denote the continuation line and tensor
$\varepsilon$ depends on the approximation errors from the tensor train
factorization and the finite difference formula.
Using the $n$-mode product introduced in
Section~\ref{subsubsec:n-mode:product}, we reformulate this
differentiation operator, in \TT{format} as:
\begin{align*}
  \left(\dfrac{\partial f}{\partial x_2}\right)^{TT}(i_1,i_2,i_3,i_4)
  = \ten{G}_1(i_1) (\ten{G}_2 \times_2 \mat{Diff})(i_2) \ten{G}_3(i_3) \ten{G}_4(i_4)
  + \varepsilon
  \quad
  \forall i_2=1,2,\ldots,(n_2-1),
\end{align*}
where matrix $\mat{Diff}$
is given by,
\begin{equation}
  \mat{Diff} \equiv
  \setlength\arraycolsep{2pt}
  \dfrac{1}{\Delta x_2}\begin{pmatrix}
    -1& 1 & & \\
    & \ddots & \ddots &\\
    & & -1 & 1
  \end{pmatrix}.
  \label{eq:diff:operator}
\end{equation}
Using this format, we apply the forward difference operation matrix
$\mat{Diff}$ only along the mode index $i_2$ of the second core
$\ten{G}_2(:,i_2,:)$ of $\ten{G}^{TT}$, see, e.g.,
Fig.~\ref{fig:TT_4D_x2_diff}.
The computational cost is greatly reduced by differentiating the multidimensional 
tensor $\ten{F}$ to differentiating only the core $\ten{G}_2$ of the \TT{tensor} 
$\ten{G}^{TT}$.
Importantly, this operation does not modify the ranks of
$\ten{G}^{TT}$; hence, no rounding operation is required to control
the rank growth.

\smallskip
\subsection{Integration Operators in \TT{format}.}
\label{APP:Int_Operators}
Consider again the function $f(x_1,x_2,x_3,x_4)$ and its full grid
tensor representation $\ten{F}(i_1,i_2,i_3,i_4)$ on a
four-dimensional, regular, Cartesian grid covering the integration
domain
$\mathcal{D}=\mathcal{D}_1\times\mathcal{D}_2\times\mathcal{D}_3\times\mathcal{D}_4$,
where each $\mathcal{D}_k$, $k=1,2,3,4$, is a 1D, bounded subinterval
of $\mathbb{R}$.
Assuming that the nodes over the domain $\mathcal{D}_2$ of $x_2$, are
chosen from a quadrature rule with corresponding weights
$\mat{w}=(w_1,w_2,\ldots,w_{n_2})$, we numerically integrate along the
independent variable $x_2$ as follows:
\begin{multline*}
  \left( \int_{\mathcal{D}_2} f(x_1,x_2,x_3,x_4) dx_2 \right)_{i_1,i_3,i_4}^{TT}
  = \sum_{\alpha_{1},\alpha_{2},\alpha_{3}=1}^{r_1,r_2,r_3} 
  \ten{G}_1(1,i_1,\alpha_1)
  \left(\sum_{i_2=1}^{n_2} w_{i_2}\ten{G}_2(\alpha_{1},i_2,\alpha_2)\right)
  \ldots\\ \ldots
  \ten{G}_3(\alpha_{2},i_3,\alpha_3)\ten{G}_4(\alpha_{3},i_4,1)
  + \varepsilon,
\end{multline*}
where the ``dots'' again denote the continuation line and tensor
$\varepsilon$ includes the approximation errors from the tensor train
factorization of $\ten{F}$ and the numerical integration.
The numerical integration along $x_2$ returns a three-dimensional
array since the index $i_2$ is absorbed by the quadrature rule
summation.
The superindex $TT$ on the left indicates that the tensor collecting
the resulting integrals still depends on indices $i_1,i_3,i_4$ and is
in the \TT{format}.
Using the contraction product introduced in
Section~\ref{subsubsec:contraction:product}, we rewrite the quadrature
rule in \TT{format} as
\begin{equation*}
  \left( \int_{x_2} f(x_1,x_2,x_3,x_4) dx_2 \right)_{i_1,i_3,i_4}^{TT}
  = \mat{G}_1(i_1)\big(\ten{G}_2\Bar{\times}_2\mat{w}\big)\mat{G}_3(i_3)\mat{G}_4(i_4).
\end{equation*}
Like numerical differentiation, numerical integration does not
modify the \TT{ranks} of $\ten{G}^{TT}$; hence, no rounding operation
is required to control rank growth.
Furthermore, the contraction product $\ten{G}_2\bar{\times}_2\mat{w}$
is an $r_1\times r_2$ matrix that can be merged to either $\ten{G}_1$
or $\ten{G}_3$ to create the three-dimensional TT-format
representation of the tensor collecting such integrals.

\subsection{Interpolation Operators in \TT{format}.}
\label{APP:Inter_Operators}
Interpolation refers to the process of estimating the values of a
function at every point that lies inside its domain of definition
using the values of that function evaluated at suitable grid nodes.
We let $\ten{I}p_\xi$ denote an interpolation operator acting only in the
$\xi$-th direction.
For example, we consider the average operator $\ten{I}p_2\ten{F}$ along
the second direction so that
\begin{align*}
  \big(\ten{I}p_2\ten{F}\big)({i_1,i_2,i_3,i_4})
  &:= \dfrac{\ten{F}(i_1,i_2+1,i_3,i_4) + \ten{F}(i_1,i_2,i_3,i_4)}{2}\\[0.5em]
  &\approx f\big(x_1(i_1), x_2(i_2)+\Delta x_2\slash{2}, x_3(i_3), x_4(i_4)\big),
\end{align*}
where $\big(x_1(i_1),x_2(i_2),x_3(i_3),x_4(i_4)\big)$ is the
coordinate vector of the grid node labeled by the multi-index
$(i_1,i_2,i_3,i_4)$ and $\Delta x_2$ is the grid stepsize along the
direction of $x_2$.
In the \TT{format}, we compute the TT-interpolation scheme 
$((\ten{I}p_2\ten{G}^{TT})\approx((\ten{I}p_2\ten{F})^{TT}$ from the values of tensor 
$\ten{G}^{TT}$ as follows:
\begin{multline*}
  \big((\ten{I}p_2\ten{G}^{TT}\big)({i_1,i_2,i_3,i_4}) :=
  \sum_{\alpha_{1},\alpha_{2},\alpha_{3}=1}^{r_1,r_2,r_3}
  \ten{G}_1(1,i_1,\alpha_1)
  \dfrac{\ten{G}_2(\alpha_{1},i_2+1,\alpha_2) +
    \ten{G}_2(\alpha_{1},i_2,\alpha_2)}{2} \dots\\ \ldots
  \ten{G}_3(\alpha_{2},i_3,\alpha_3)\ten{G}_4(\alpha_{3},i_4,1) +
  \varepsilon,
\end{multline*}
where the ``dots'' again denote the continuation line and
$\varepsilon$ is the error depending on the \TT{factorization} and the
interpolation scheme.
The $n$-mode product of Section~\eqref{subsubsec:n-mode:product} makes it 
possible to reformulate the action of the interpolation operator
$(\ten{I}p_2$ in the compact \TT{matrix format} as:
\begin{equation*}
  \big((\ten{I}p_2\ten{G}^{TT}\big)(i_1,i_2,i_3,i_4) =
  \mat{G}_1(i_1)\big(\ten{G}_2\times_2\mat{I_2}\big)(i_2)\mat{G}_3(i_3)\mat{G}_4(i_4),
\end{equation*}
where
\begin{equation*}
  \mat{I_p} \equiv
  \setlength\arraycolsep{2pt}
  \dfrac{1}{2}\begin{pmatrix}
    1& 1 & & \\
    & \ddots & \ddots &\\
    & & 1 & 1
  \end{pmatrix}.
\end{equation*}
Again, the discrete separation of the mode indices provided by the
TT-format representation allows us to apply the interpolation matrix
$\mat{I_2}$ to the mode index $i_2$ of the second core
$\ten{G}_2(:,i_2,:)$.
As noted in the case of numerical differentiation and integration,
tensor $(\ten{I}p_2\ten{F}^{TT}$ is already in \TT{format} and the
interpolation operation does not modify the ranks of $\ten{G}^{TT}$;
hence, no rounding operation is required to reduce the \TT{ranks}.

\end{document}